\documentclass{elsart}
\usepackage{amsmath}
\numberwithin{equation}{section}

\newtheorem{dl}{Theorem}
\newtheorem{yl}{Lemma}
\newtheorem{dy}{Definition}

\newtheorem{tl}{Corollary}
\newtheorem{lz}{Example}

\theoremstyle{definition}

\theoremstyle{remark}
\begin{document}
\begin{frontmatter}
\title{The $(f,g)$-inversion formula and its applications: the $(f,g)$-summation formula}
\author{Xinrong Ma}
\address[P.R.China]{Department of Mathematics, SuZhou University, SuZhou 215006, P.R.China}
\begin{abstract}
A complete characterization of two functions $f(x,y)$ and $g(x,y)$
in the $(f,g)$-inversion is presented. As an application to the
theory of hypergeometric series, a general bibasic summation formula
determined by $f(x,y)$ and $g(x,y)$ as well as four arbitrary
sequences is obtained which unifies Gasper and Rahman's, Chu's and
Macdonald's bibasic summation formula. Furthermore, an alternative
proof of the $(f,g)$-inversion derived from the $(f,g)$-summation
formula is presented. A bilateral $(f,g)$-inversion containing
Schlosser's bilateral matrix inversion as a special case  is also
obtained.
\end{abstract}
\begin{keyword}
 Bilateral formal power series, $(f,g)$-inversion, $(f,g)$-expansion
formula, basic, elliptic, bibasic, bilateral hypergeometric series,
transformation formula, $(f,g)$-summation formula.
\end{keyword}
\end{frontmatter}
\maketitle

\maketitle
\section{Introduction}
\setcounter{equation}{0} \setcounter{dl}{0}
\setcounter{tl}{0}\setcounter{yl}{0}\setcounter{lz}{0} As
well-known, matrix inversions over the complex field $\mathcal{C}$,
 play a very important
role in the theory of basic hypergeometric series. Recall that
 a matrix inversion is define to be a pair of infinite lower triangular matrices
$F=(f_{n,k})_{n,k \in Z}$ and $G=(g_{n,k})_{n,k \in Z}$
 satisfying
 $$\displaystyle \sum_{i=k}^nf_{n,i}g_{i,k}=\delta_{n,k},
$$
where $\delta$ denotes the Kronecker delta, $Z$ denotes the set of
integers.
 If such a pair of matrices $F$ and $G$ are not lower triangular, we call it bilateral matrix inversion.
 In our
previous paper \cite{0020}, we have established the following matrix
inversion, named the $(f,g)$-inversion.
 \begin{dl} \label{math4} Let $f(x,y)$ and $g(x,y)$ be two arbitrary functions over the complex field $\mathcal{C}$ in variables
 $x,y$. Suppose that $g(x,y)$ is antisymmetric, i.e., $g(x,y)=-g(y,x)$.
 Let $F=(f_{n,k})_{n,k\in Z}$ and $G=(g_{n,k})_{n,k\in Z}$ be two matrices
with entries given by
\begin{eqnarray}
f_{n,k}&=&\frac{\prod_{i=k}^{n-1}f(x_i,b_k)}
{\prod_{i=k+1}^{n}g(b_i,b_k)}\label{a12}\qquad\mbox{and}\\
g_{n,k}&=&
\frac{f(x_k,b_k)}{f(x_n,b_n)}\frac{\prod_{i=k+1}^{n}f(x_i,b_n)}
{\prod_{i=k}^{n-1}g(b_i,b_n)},\label{b12}\quad\mbox{respectively}.
\end{eqnarray}
 Then $F=(f_{n,k})_{n,k\in Z}$ and $G=(g_{n,k})_{n,k\in Z}$ is a
matrix inversion if and only if $f(x,y)\in \mbox{\sl
Ker}\mathcal{L}^{(g)}_{3}$ or $f(x,y)\in \mbox{\sl
Ker}\mathcal{L}_{3}$, where $\mbox{\sl Ker}\mathcal{L}_{3}^{(g)}$
and $\mbox{\sl Ker}\mathcal{L}_{3}$ denote the sets of functions
 $f(x,y)$ such that for all $ a,b,c,x\in \mathcal{C}$
\begin{eqnarray}
&&\label{add}g(a,b)f(x,c)-g(a,c)f(x,b)+g(b,c)f(x,a)=0\quad\mbox{and}\\
&&f(a,b)f(x,c)-f(a,c)f(x,b)+f(b,c)f(x,a)=0,\label{add0}
\end{eqnarray}
 respectively.
\end{dl}
Our results in \cite{0020} state that the  $(f,g)$-inversion covers
 all known matrix inversions such as the Gould-Hsu formula, Krattenthaler's  inversion formula,
 as well as Warnaar's elliptic
 matrix inversion.
The aim of this paper, as the further developments in work of our
paper \cite{0020}, is twofold. First, we will give a complete
characterization of $f$ and $g$ such that $f\in \mbox{\sl
Ker}\mathcal{L}_{3}^{(g)}$, which was left as an open problem in
\cite{0020}. The answer is presented in Section 2. Second, we will
set up a rather general summation formula via such a pair of
functions $f$ and $g$. Later as we will see, this summation formula
unifies all known bibasic summation formulas which are extensions of
Jackson's q-analogue of Whipple's summation formula for a
terminating very-well-poised balanced $\,_8\phi_7$ \cite[2.6.2]{10}
and are due to several authors
(cf.\cite{alver,macd,chu,10,101,new1,1000}). For a good survey about
this kind of summation formulae, we refer the reader to Gasper and
Rahman's book \cite[Sections 3.6-3.7]{10}
 and for their $U(n+1)$-extensions to \cite{macd,milne1}. In
addition, from the $(f,g)$-summation formula, we obtain a very short
proof of the $(f,g)$-inversion and a bilateral $(f,g)$-inversion
containing Schlosser's bilateral matrix inversion as a special case,
both are presented in Section 3.

 {\sl Notations and conventions.} Throughout this
paper, for convenience, we will mainly work on the setting of formal
bilateral series over the complex field $C$, and use the notation
 $\mathcal{C}[x,y]$ to denote the ring of formal bilateral series
in $x,y$. As usual, we employ the following standard notations for
the q-shifted factorial: \begin{eqnarray*} &&
(a;q)_n=\prod_{i=0}^{n-1}(1-aq^i),
(a;q)_{\infty}=\prod_{i=0}^{\infty}(1-aq^i)\\
&&
 (a_1,a_2,\cdots,a_m;q)_n=(a_1;q)_n(a_2;q)_n\cdots(a_m;q)_n.
   \end{eqnarray*}

\section{The complete characterizations of $\mbox{\sl Ker}\mathcal{L}_{3}$ and $\mbox{\sl Ker}\mathcal{L}^{(g)}_{3}$}
In this section, we will characterize two nonzero functions $f(x,y)$
and $g(x,y)$ required by the $(f,g)$-inversion. For this, we first
introduce some notation and preliminaries.
\begin{dy}\label{dy1}Let $f(x,y)$ and $g(x,y)$ be two arbitrary nonzero functions over $\mathcal{C}$ in variables
 $x,y$. Suppose that for all four numbers  $a,b,c,x$,
\begin{eqnarray*}
g(a,b)f(x,c)-g(a,c)f(x,b)+g(b,c)f(x,a)=0,
\end{eqnarray*}
then $f(x,y)$ is called orthogonal to $g(x,y).$ Further, $f(x,y)$ is
called self-orthogonal if $f(x,y)$ is orthogonal to itself.
\end{dy}
Write $f(x,y)\bot g(x,y)$, or in short, $f\bot g$ if $f(x,y)$ is
orthogonal to $g(x,y)$. As it stands,  $f\bot g$ does not mean
$g\bot f$. This definition allows us to rewrite
$$
\mbox{\sl Ker}\mathcal{L}^{(g)}_{3}=\{f|f\bot g\}; \mbox{\sl
Ker}\mathcal{L}_{3}=\{f|f\bot f\}.
$$

Since the set $ \mbox{\sl Ker}\mathcal{L}^{(g)}_{3}$ is of value to
the $(f,g)$-inversion, it seems necessary  to find the explicit
expressions of $f$ and $g$. The possibility for us to do this is
within the ring of formal bilateral series. At first, we need
\begin{yl}\label{2.1}Let
$f(x,y)=\sum_{i,j=-\infty}^{\infty}\lambda(i,j)x^iy^j\in C[x,y]$.
Then $f\bot f$ if and only if  for any
 integers $ m,i,j,k\in Z,$
\begin{eqnarray}
\lambda(m,i)\lambda(k,j)-\lambda(k,i)\lambda(m,j)+
\lambda(k,m)\lambda(i,j)=0.\label{xiugai}
\end{eqnarray}
\end{yl}
Based on (\ref{xiugai}), it is easy to show that $f(x,y)$ is
antisymmetric, i.e., $f(x,y)=-f(y,x)$.

\begin{yl}Let $g(x,y)=\sum_{i,j=-\infty}^{\infty}c(i,j)x^iy^j,f(x,y)=\sum_{i,j=-\infty}^{\infty}\lambda(i,j)x^iy^j$. Then
$f\bot g$ if and only if  for any
 integers $ m,i,j,k\in Z$,
\begin{eqnarray}
  c(m,i)\lambda(k,j)-c(m,j)\lambda(k,i)+
c(i,j)\lambda(k,m)=0.
\end{eqnarray}
\end{yl}

 Before we  give the explicit expressions of $f$ and $g$, it is better for us to work out whether there
exists any connection between $g\bot g$ and $f\bot g$. This idea
allows us to set up
\begin{yl}\label{6.4} Let
$g(x,y)=\sum_{i,j=-\infty}^{\infty}c(i,j)x^iy^j\neq 0$. Then $g\bot
g$ if and only if there are at least two integers $m_0,k_0$, such
that $c(m_0,k_0)\neq 0$ and such that for arbitrary $i,j\in Z$,
there holds that
\begin{eqnarray}
&&c(m_0,k_0)c(i,j)-c(m_0,i)c(k_0,j)+ c(m_0,j)c(k_0,i)=0.\label{0000}
\end{eqnarray}
\end{yl}
{\sl Proof.} Obviously, since $g\neq 0$, there is at least one
nonzero coefficient $c(m_0,k_0)$. Let $g\bot g$. Then by Definition
 \ref{dy1}, we see that (\ref{0000}) holds. Now, assume that (\ref{0000})
holds for arbitrary integers $i,j$. Without any loss of generality,
let $c(m_0,k_0)=1$. Thus, for any quadruple of fixed $i,j,i_1,j_1\in
Z$, there must hold that
\begin{eqnarray*}\left\{%
\begin{array}{ll}
&c(i,j)=c(m_0,j)c(i,k_0)-c(k_0,j)c(i,m_0); \\
 &c(i_1,j_1)=c(m_0,j_1)c(i_1,k_0)-c(k_0,j_1)c(i_1,m_0); \\
 &c(i,i_1)=c(m_0,i_1)c(i,k_0)-c(k_0,i_1)c(i,m_0); \\
  &c(j,j_1)=c(m_0,j_1)c(j,k_0)-c(k_0,j_1)c(j,m_0); \\
 & c(i,j_1)=c(m_0,j_1)c(i,k_0)-c(k_0,j_1)c(i,m_0); \\
  &c(i_1,j)=c(m_0,j)c(i_1,k_0)-c(k_0,j)c(i_1,m_0).
  \end{array}%
\right.
\end{eqnarray*}
Set
\begin{eqnarray*} &&c(m_0,i)=x_1,  c(m_0,j)=x_2, c(m_0,i_1)=x_3,
c(m_0,j_1)=x_4;\\
&& c(k_0,i)=-y_1, c(k_0,j)=-y_2, c(k_0,i_1)=-y_3, c(k_0,j_1)=-y_4.
\end{eqnarray*}
Insert these into the preceding relations to arrive at
\begin{eqnarray*}&&c(i,j)=x_2y_1-x_1y_2, c(i,i_1)=x_3y_1-x_1y_3, c(i,j_1)=x_4y_1-x_1y_4;
\\
&& c(i_1,j_1)=x_4y_3-x_3y_4, c(j,j_1)=x_4y_2-x_2y_4,
c(i_1,j)=x_2y_3-x_3y_2.
     \end{eqnarray*}
And then by a bit of straightforward calculation, it is easily seen
that
\begin{eqnarray*}
  c(i_1,j_1)c(i,j)-c(i_1,j)c(i,j_1)+c(j_1,j)c(i,i_1)=0,
\end{eqnarray*}
which gives that $g\bot g$.

The next lemma gives the relationship between $g\bot g$ and $ f\bot
g$.
 \begin{yl}\label{1.2}
Let $f(x,y)=\sum_{i,j=-\infty}^{\infty}\lambda(i,j)x^iy^j,
g(x,y)=\sum_{i,j=-\infty}^{\infty}c(i,j)x^iy^j\in C[x,y]$  be two
nonzero
 bilateral series, $g(x,y)=-g(y,x)$. Then  $g\bot g$ provided that $ f\bot g$.
\end{yl}
{\sl Proof.} Since $g(x,y)$ is a nonzero bilateral series in $x,y$,
 there exists  at least one nonzero coefficient $c(m_0,k_0)$. Under this condition, solving  $f\bot g$ , i.e.,
\begin{eqnarray}c(m,k)\lambda(i,j)-
c(m,j)\lambda(i,k)+c(k,j)\lambda(i,m)=0\label{qianshi}
\end{eqnarray}
 immediately gives that for arbitrary $i,j$,\begin{eqnarray}\lambda(i,j)=
\frac{c(m_0,j)\lambda(i,k_0)+c(j,k_0)\lambda(i,m_0)}{c(m_0,k_0)}.\label{biaoshi}
\end{eqnarray}
Substitute each $\lambda(i,j)$  in (\ref{qianshi}) by
(\ref{biaoshi}) and then rearrange the corresponding result and find
\begin{eqnarray*}
&&\frac{\lambda(k,k_0)}{c(m_0,k_0)}
\left[c(m,i)c(m_0,j)-c(m,j)c(m_0,i)+c(i,j)c(m_0,m)\right]
\\&+&
   \frac{\lambda(k,m_0)}{c(m_0,k_0)}\left[c(m,i)c(j,k_0)-c(m,j)c(i,k_0)+c(i,j)c(m,k_0)\right]=0.
\end{eqnarray*}
Note that $k,m,i,j$ are arbitrary and $f(x,y)$ is nonzero, there are
at least one integer $k$, such that $\lambda(k,k_0)\neq 0$ or
$\lambda(k,m_0)\neq 0$. Without any loss of generality, suppose
$\lambda(k,k_0)\neq 0$. Now, set $j=k_0$. Since $c(k_0,k_0)=0$,
which follows from $g(x,y)=-g(y,x)$, it is easily seen that the
second sum in the left-hand side of this identity vanishes. So
 the result is that for arbitrary $i,m\in Z$,
\begin{eqnarray*}
\frac{\lambda(k,k_0)}{c(m_0,k_0)}
\left[c(m,i)c(m_0,k_0)-c(m,k_0)c(m_0,i)+c(i,k_0)c(m_0,m)\right]=0,
\end{eqnarray*}
 which is equivalent to
\begin{eqnarray*}
c(m,i)c(m_0,k_0)-c(m,k_0)c(m_0,i)+c(i,k_0)c(m_0,m)=0.
\end{eqnarray*}
 By Lemma \ref{6.4}, it is easily seen that  $g(x,y)\in
\mbox{\sl Ker}\mathcal{L}_3$.

 Lemma \ref{1.2} provides a necessary condition for two functions $f$ and $g$ such that $f
 \bot g.$
These lemmas in together  makes it possible to find the explicit
expressions for such $f$ or $g$. In what follows,
 it is always assumed that $g\bot g$.

\begin{dl}\label{tezhdl1}Suppose that
$f(x,y)=\sum_{i,j=-\infty}^{+\infty}\lambda(i,j)x^iy^j\in C[x,y]$.
Then $f\bot f$ if and only if there exist two complex sequences
$\{p_i\}_{i\in Z}$ and $\{q_i\}_{i\in Z}$ such that
$p_{m_0}=0,q_{k_0}=0$, and $p_{k_0}=-q_{m_0}\neq 0$,
\begin{eqnarray}
f(x,y)=P(x)Q(y)-P(y)Q(x),
\end{eqnarray}
where
$P(x)=\sum_{i=-\infty}^{+\infty}p_ix^i,Q(x)=\frac{1}{q_{m_0}}\sum_{i=-\infty}^{+\infty}q_ix^i$.
\end{dl}
{\sl Proof.} By Lemma \ref{2.1}, it is easily seen that $f\bot f$ is
equivalent to
\begin{eqnarray}
\lambda(i,j)\lambda(m,k)+
\lambda(i,k)\lambda(j,m)-\lambda(i,m)\lambda(j,k)=0.\label{4.2}
\end{eqnarray}
Assume that  $\lambda(m_0,k_0)\neq 0$. In this case, it is easy to
solve (\ref{4.2}) and to get
\begin{eqnarray}\lambda(i,j)&=&\frac{\lambda(i,m_0)\lambda(j,k_0)-
\lambda(i,k_0)\lambda(j,m_0)}{\lambda(m_0,k_0)}.\label{4.4}\end{eqnarray}
Starting with  this identity, we proceed to show the theorem in
question by proving  that (\ref{4.2}) is  equivalent to (\ref{4.4})
and then expressing the solution given by (\ref{4.4}) in terms of
generating function. At first, by the definition,
 it is easily seen that  $(\ref{4.2})$ implies $(\ref{4.4})$. We need only  to show is
 that  $(\ref{4.2})$ can be reduced from  $(\ref{4.4})$. To do this,
 replacing each $\lambda(\bullet,\bullet)$ by (\ref{4.4}) in the left-hand side of
  (\ref{4.2}) to calculate
\begin{eqnarray*}\lefteqn{\lambda(i,j)\lambda(m,k)+
\lambda(i,k)\lambda(j,m)-\lambda(i,m)\lambda(j,k)}\\&=&
\frac{\lambda(i,m_0)\lambda(j,k_0)-\lambda(i,k_0)\lambda(j,m_0)}{\lambda(m_0,k_0)}\cdot
\frac{\lambda(m,m_0)\lambda(k,k_0)-\lambda(m,k_0)\lambda(k,m_0)}{\lambda(m_0,k_0)}\\&-&
\frac{\lambda(i,m_0)\lambda(k,k_0)-\lambda(i,k_0)\lambda(k,m_0)}{\lambda(m_0,k_0)}\cdot
\frac{\lambda(j,m_0)\lambda(m,k_0)-\lambda(j,k_0)\lambda(m,m_0)}{\lambda(m_0,k_0)}\\&+&
\frac{\lambda(i,m_0)\lambda(m,k_0)-\lambda(i,k_0)\lambda(m,m_0)}{\lambda(m_0,k_0)}\cdot
\frac{\lambda(j,m_0)\lambda(k,k_0)-\lambda(j,k_0)\lambda(k,m_0)}{\lambda(m_0,k_0)}
.
\end{eqnarray*}
Now, expand and rearrange the right-hand side of this identity to
arrive
 at
  \begin{eqnarray*}\lefteqn{\lambda(i,j)\lambda(m,k)+
\lambda(i,k)\lambda(j,m)-\lambda(i,m)\lambda(j,k)}\\&=&
\frac{\lambda(i,m_0)\lambda(k,k_0)}{\lambda(m_0,k_0)^2}
[\lambda(j,k_0)\lambda(m,m_0)+\lambda(j,m_0)\lambda(m,k_0)-\lambda(m,k_0)\lambda(j,m_0)]\\&-&
\frac{\lambda(i,m_0)\lambda(j,k_0)}{\lambda(m_0,k_0)^2}
[\lambda(m,k_0)\lambda(k,m_0)+\lambda(k,k_0)\lambda(m,m_0)-\lambda(m,k_0)\lambda(k,m_0)]\\&-&
\frac{\lambda(i,k_0)\lambda(j,m_0)}{\lambda(m_0,k_0)^2}
[\lambda(m,m_0)\lambda(k,k_0)+\lambda(k,m_0)\lambda(m,k_0)-\lambda(m,m_0)\lambda(k,k_0)]\\&+&
\frac{\lambda(i,k_0)\lambda(k,m_0)}{\lambda(m_0,k_0)^2}
[\lambda(j,m_0)\lambda(m,k_0)+\lambda(j,k_0)\lambda(m,m_0)-\lambda(m,m_0)\lambda(j,k_0)].
\end{eqnarray*}
 Obviously, each term within the curly
braces can be further simplified. The corresponding result is
\begin{eqnarray*}\lefteqn{\lambda(i,j)\lambda(m,k)+
\lambda(i,k)\lambda(j,m)-\lambda(i,m)\lambda(j,k)}\\&=&
\frac{\lambda(i,m_0)\lambda(k,k_0)\lambda(j,k_0)\lambda(m,m_0)}{\lambda(m_0,k_0)^2}-
   \frac{\lambda(i,m_0)\lambda(j,k_0)\lambda(k,k_0)\lambda(m,m_0)}{\lambda(m_0,k_0)^2}-\\& &
   \frac{\lambda(i,k_0)\lambda(j,m_0)\lambda(k,m_0)\lambda(m,k_0)}{\lambda(m_0,k_0)^2}+
   \frac{\lambda(i,k_0)\lambda(k,m_0)\lambda(j,m_0)\lambda(m,k_0)}{\lambda(m_0,k_0)^2}\\
&=&0.
\end{eqnarray*}
 Finally, we turn to express
$\lambda(i,j)$ given by (\ref{4.4}) in terms of generating function.
To do this, suppose that $\{p_i\}_{i\in Z}$ and $\{q_i\}_{i\in Z}$
are
 arbitrary complex sequences such that $p_{k_0}=-q_{m_0}\neq 0$. Define that
\begin{eqnarray*}\lambda(i,m_0)=p_i,\lambda(i,k_0)=q_i,i\in Z.\end{eqnarray*}
Note that $p_{m_0}=0$,\ $q_{k_0}=0$, because $f(x,y)=-f(y,x)$.
Hence, $\lambda(i,j)$ given by (\ref{4.4}) may be rewritten as
\begin{eqnarray}f(x,y)=P(x)Q(y)-P(y)Q(x),
\end{eqnarray}
where
$P(x)=\sum_{i=-\infty}^{+\infty}p_ix^i,Q(x)=\frac{1}{q_{m_0}}\sum_{i=-\infty}^{+\infty}q_ix^i$.
This gives the complete proof of theorem.

\begin{tl}\label{coll21} The following functions $f(x,y)$ are self orthogonal.
\begin{eqnarray*}
f(x,y)=x-y,(y-x)(1-\frac{xy}{d}),(x-y)(1-\frac{b}{axy}).
\end{eqnarray*}
\end{tl}
As an important case that both $f(x,y)\in \mbox{\sl
Ker}\mathcal{L}_{3}$ and $f(x,y)$ is infinite bilateral series, we
obtain an alternate proof of the {\sl addition formula} of the theta
function.
\begin{tl} \label{coll22} Define
$\theta(x)=(x;q)_\infty(\frac{q}{x};q)_\infty\,\, (|q|<1)$, and
$f(x,y)=y\theta(xy)\theta(\frac{x}{y})$. Then $f(x,y)$ is self
orthogonal.
\end{tl}
{\sl Proof.} In Theorem \ref{tezhdl1}, set $(m_0,k_0)=(1,0)$ and
\begin{eqnarray*}p_i&=&
\left\{\begin{array}{lll} &\displaystyle\frac{q^{m^2-m}
}{(q;q)_\infty^2},&\mbox{if}\,\, i =2m;\\
&0,& \mbox{if}\,\,i=2m+1
\end{array}\right.\quad \mbox{and}\quad q_i=\left\{\begin{array}{lll} &\displaystyle- \frac{q^{m^2}
}{(q;q)_\infty^2},&\mbox{if}\,\, i=2m+1;\\
&0,&\mbox{if}\,\, i=2m.
\end{array}\right.
\end{eqnarray*}
A direct calculation gives the solution
\begin{eqnarray*}f(x,y)&=&\frac{1}{(q;q)_\infty^2}\sum_{i,j=-\infty}^{+\infty}
 q^{i^2-i+j^2}(x^{2i}y^{2j+1}-y^{2i}x^{2j+1}).\end{eqnarray*}
By the definition of $\theta(x)$ and Jacobi's \emph{triple product
identity}
\begin{eqnarray*}
\sum_{i=-\infty}^{+\infty}(-1)^iq^{i \choose
2}x^i=(x;q)_\infty(\frac{q}{x};q)_\infty(q;q)_\infty,
\end{eqnarray*}
it follows that $f(x,y)=y\theta(xy)\theta(\frac{x}{y}).$

 The next theorem gives the explicit expression of $f(x,y)\in \mbox{\sl
Ker}\mathcal{L}^{(g)}_{3}$.
\begin{dl}\label{tezhdl2}
Let $g(x,y)=\sum_{i,j=-\infty}^{\infty}c(i,j)x^iy^j\neq 0$ such that
$g\bot g$. Then $ f\bot g$ if and only if there exists
$c(m_0,k_0)\neq 0$, such that
\begin{eqnarray}
f(x,y)&=&\frac{1}{c(m_0,k_0)}\left(P(x)[x^{m_0}]g(x,y)-Q(x)[x^{k_0}]g(x,y)\right),\label{new}
\end{eqnarray}
where $P(x)=\sum_{i=-\infty}^{\infty}p_ix^i,\,
Q(x)=\sum_{i=-\infty}^{\infty}q_ix^i,$ $p_i,q_i$ are arbitrary
sequences, $[x^{m_0}]g(x,y)$ stands the coefficient of $x^{m_0}$ in
$g(x,y)$.
\end{dl}
{\sl Proof.}
 Suppose that
 $f(x,y)=\sum_{i,j=-\infty}^{\infty}\lambda(i,j)x^iy^j.$
Observe that $f\bot g$ is equivalent to, by Lemma \ref{1.2}, the
relation
\begin{eqnarray}
\left\{%
\begin{array}{ll}
     &c(m,i)c(k,j)-c(k,i)c(m,j)+c(k,m)c(i,j)=0;\\
         &c(m,i)\lambda(k,j)-
c(m,j)\lambda(k,i)+c(i,j)\lambda(k,m)=0.
\end{array}%
\right.    \label{233}
\end{eqnarray}
 Let
$m_0$ and $k_0$ be two integers such that $c(m_0,k_0)\neq 0$. Use
this to solve (\ref{233}) and  find
that\begin{eqnarray}\lambda(k,j)=
\frac{c(m_0,j)\lambda(k,k_0)+c(j,k_0)\lambda(k,m_0)}{c(m_0,k_0)}.\label{2.4}
\end{eqnarray}
Proceeding as before, we need only to show that (\ref{233}) and
(\ref{2.4}) are equivalent to each other. At first, $(\ref{233})$
implies $(\ref{2.4})$, which is obvious. Conversely, $(\ref{2.4})$
implies $(\ref{233})$, too. To show this, replacing each
$\lambda(\bullet,\bullet)$  by (\ref{2.4}) in the left-hand side of
the second identity of (\ref{233}) directly gives that
\begin{eqnarray*}\lefteqn{c(m,i)\lambda(k,j)-
c(m,j)\lambda(k,i)+c(i,j)\lambda(k,m)}\\
&=&c(m,i)\frac{c(m_0,j)\lambda(k,k_0)+c(j,k_0)\lambda(k,m_0)}{c(m_0,k_0)}
\\&-& c(m,j)\frac{c(m_0,i)\lambda(k,k_0)+c(i,k_0)\lambda(k,m_0)}{c(m_0,k_0)}\\
&+&c(i,j)\frac{c(m_0,m)\lambda(k,k_0)+c(m,k_0)\lambda(k,m_0)}{c(m_0,k_0)}\\
&=&\frac{\lambda(k,k_0)}{c(m_0,k_0)}
\left[c(m,i)c(m_0,j)-c(m,j)c(m_0,i)+c(i,j)c(m_0,m)\right]
\\&+&
   \frac{\lambda(k,m_0)}{c(m_0,k_0)}\left[c(m,i)c(j,k_0)-c(m,j)c(i,k_0)+c(i,j)c(m,k_0)\right].
\end{eqnarray*}
Observe that $g(x,y)\in \mbox{\sl Ker}\mathcal{L}_3$ gives
\begin{eqnarray*}\left\{%
\begin{array}{ll}
&c(m,i)c(m_0,j)-c(m,j)c(m_0,i)+c(i,j)c(m_0,m)=0;\\
&c(m,i)c(j,k_0)-c(m,j)c(i,k_0)+c(i,j)c(m,k_0)=0.
\end{array}
\right.
\end{eqnarray*}
Thus, the previous identity becomes
\begin{eqnarray*}c(m,i)\lambda(k,j)-c(m,j)\lambda(k,i)+c(i,j)\lambda(k,m)=0.
\end{eqnarray*}
So (\ref{233}) follows.  Therefore, (\ref{2.4}) gives the desired
function. Next, let $\{p_i\}$ and $\{q_i\}$ be arbitrary
 complex sequences. Define
 $$P(x)=\sum_{i=-\infty}^{\infty}p_ix^i,Q(x)=\sum_{i=-\infty}^{\infty}q_ix^i,$$ and
 $\lambda(i,k_0)=p_i,\lambda(j,m_0)=q_j.$
Then (\ref{2.4}) can be rewritten as  \begin{eqnarray*}
\lambda(i,j)=\frac{c(m_0,j)p_i+c(j,k_0)q_i}{c(m_0,k_0)},
\end{eqnarray*}
which leads to \begin{align}
f(x,y)=\frac{1}{c(m_0,k_0)}\sum_{i,j=-\infty}^{\infty}
(c(m_0,j)p_i+c(j,k_0)q_{i})x^iy^j\nonumber\\
=\frac{1}{c(m_0,k_0)}\sum_{i=-\infty}^{\infty}
p_ix^i\sum_{j=-\infty}^{\infty} c(m_0,j)y^j-\frac{1}{c(m_0,k_0)}
\sum_{i=-\infty}^{\infty}q_{i}x^i\sum_{j=-\infty}^{\infty}
c(k_0,j)y^j\nonumber\\
=\frac{1}{c(m_0,k_0)}\left(P(x)[x^{m_0}]g(x,y)-Q(x)[x^{k_0}]g(x,y)\right).
\end{align}

Note that the above results in together provides us with a
constructive way to derive different $(f,g)$-inversions and
expansion formulae, when this is necessary.

\begin{tl}\label{coll23} The following function pairs $(f,g)$ satisfy that $f\in \mbox{\sl
Ker}\mathcal{L}^{(g)}_{3}$.
\begin{eqnarray*}
(f,g)&=&(P(x)+y Q(x),x-y),\left((1-axy)(1-b\frac{x}{y}),(x-y)(1-\frac{b}{axy})\right),\\
&&\left((x+y)(x+\frac{b}
  {ay}),(x-y)(1-\frac{b}{axy})\right).
\end{eqnarray*}
\end{tl}
{\sl Proof.} By Corollary \ref{coll21}, each $g\in \mbox{\sl
Ker}\mathcal{L}_{3}$. It only needs to verify $f\bot g$ by
\emph{Maple} in a straightforward way. Left to the reader.
\section{The $(f,g)$-summation
formula} As mentioned in Section 1, we start with a general
summation formula related with $f$ and $g$ satisfying $f\bot g$.
\begin{yl}\label{13yl} Let $\{a_i\},\{b_i\},\{c_i\},\{d_i\}$ be arbitrary
sequences such that none of the denominators in (\ref{news1155})
vanish. Then for any integer $m\geq 0$,
\begin{align}
\sum^m_{k=0} f(a_k,b_k)g(c_k,d_k)\frac{\prod_{j=1}^{k-1}f(a_j,c_j)}
{\prod_{j=1}^{k} f(a_j,d_j)}\frac{\prod_{j=1}^{k-1}g(b_j,d_j)}
{\prod_{j=1}^{k}g(b_j,c_j)}\nonumber\\=\frac{\prod_{j=1}^{m}f(a_j,c_j)}
{\prod_{j=1}^{m}f(a_j,d_j)}\frac{\prod_{j=1}^{m} g(b_j,d_j)}
{\prod_{j=1}^{m}g(b_j,c_j)}+C\label{news1155}
\end{align}
provided $f\bot g$, i.e., $f\in\mbox{\sl
Ker}\mathcal{L}^{(g)}_{3}$,where the constant $C$ is uniquely
decided by the initial conditions.
\end{yl}
The idea of the proof relies on the difference equation: let $\Delta
S_n=S_{n+1}-S_n =t_n$. Then by telescoping, $\sum_{k=0}^nt_k=S_n+C$.
Conversely, if $\sum_{k=0}^nt_k=S_n+C$ and $t_k$ ($k\neq n)$ is
independent of $n$, then $t_n=\Delta\,S_{n}$. The constant $C$ is
given by the initial condition.

 {\sl Proof.}
Define
\begin{eqnarray}
 F(m)&=&\frac{\prod_{j=1}^{m}f(a_j,c_j)}
{\prod_{j=1}^{m}f(a_j,d_j)}\frac{\prod_{j=1}^{m} g(b_j,d_j)}
{\prod_{j=1}^{m}g(b_j,c_j)},\end{eqnarray} and assume further that
\begin{eqnarray}
F(m)&=&\sum^m_{k=0} f(a_k,b_k)G(k)\frac{\prod_{j=0}^{k-1}g(b_i,d_j)}
{\prod_{j=1}^{k} f(x_i,d_j)}+C.\label{newh}
 \end{eqnarray}
 Thus, as mentioned above, we need only to show that $G(n)$ is independent of
$m$. We proceed to do this in a straightforward calculation. Suppose
\begin{eqnarray}
\label{111}\frac{\prod_{j=1}^{m}f(a_j,c_j)}
{\prod_{j=1}^{m}f(a_j,d_j)}\frac{\prod_{j=1}^{m} g(b_j,d_j)}
{\prod_{j=1}^{m}g(b_j,c_j)}=\sum^m_{k=0}
f(a_k,b_k)G(k)\frac{\prod_{j=0}^{k-1}g(b_i,d_j)} {\prod_{j=1}^{k}
f(x_i,d_j)}+C.\end{eqnarray}
 Calculation with the assumption  gives that
\begin{eqnarray}G(k)&=&\frac{F(k)-F(k-1)}{f(a_k,b_k)}
\frac{\prod_{j=1}^{k}f(a_j,d_j)}{\prod_{j=0}^{k-1}g(b_j,d_j)}\nonumber
\\ &=&\frac{1}{f(a_k,b_k)}
\left\{\frac{F(k)}{F(k-1)}-1\right\}F(k-1)
\frac{\prod_{j=1}^{k}f(a_j,d_j)}{\prod_{j=0}^{k-1}g(b_j,d_j)}.\label{expandefcoeffff}
\end{eqnarray}
Observe that the term within the curly braces
\begin{eqnarray*}
  \frac{F(k)}{F(k-1)}-1=\frac{f(a_k,c_k)}
{g(b_k,c_k)}\frac{ g(b_k,d_k)} {f(a_k,d_k)}-1
=\frac{f(a_k,c_k)g(b_k,d_k)-
g(b_k,c_k)f(a_k,d_k)}{g(b_k,c_k)f(a_k,d_k)}.
\end{eqnarray*}
Simplify the enumerator in the latter expression further by means of
 the known condition that $f\bot g$, namely say,
\begin{equation}
\label{use} f(a_k,c_k)g(b_k,d_k)-
g(b_k,c_k)f(a_k,d_k)=f(a_k,b_k)g(c_k,d_k),
\end{equation}
to get
\begin{eqnarray*}
  \frac{F(k)}{F(k-1)}-1=\frac{f(a_k,b_k)g(c_k,d_k)}{g(b_k,c_k)f(a_k,d_k)}.
\end{eqnarray*}
 The production after some calculation is
\begin{eqnarray*}
F(k-1) \frac{\prod_{j=1}^{k}f(a_j,d_j)}{\prod_{j=0}^{k-1}
g(b_j,d_j)}&=&\frac{\prod_{j=1}^{k-1}f(a_j,c_j)}
{\prod_{j=1}^{k-1}g(b_j,c_j)}\frac{\prod_{j=1}^{k-1} g(b_j,d_j)}
{\prod_{j=1}^{k-1}f(a_j,d_j)}\frac{\prod_{j=1}^{k}f(a_j,d_j)}{\prod_{j=0}^{k-1}
g(b_j,d_j)}\\
&=&\frac{f(a_k,d_k)}{g(b_0,d_0)}\frac{\prod_{j=1}^{k-1}f(a_j,c_j)}
{\prod_{j=1}^{k-1}g(b_j,c_j)}.
\end{eqnarray*}
Insert all these notes into (\ref{expandefcoeffff}). It gives
\begin{eqnarray*}G(k)&=&\frac{1}{f(a_k,b_k)}
\frac{f(a_k,b_k)g(c_k,d_k)}{g(b_k,c_k)f(a_k,d_k)}
\frac{f(a_k,d_k)}{g(b_0,d_0)}\frac{\prod_{j=1}^{k-1} f(a_j,c_j)}
{\prod_{j=1}^{k-1}g(b_j,c_j)}\\
&=& \frac{g(c_k,d_k)}{g(b_0,d_0)}\frac{\prod_{j=1}^{k-1}f(a_j,c_j)}
{\prod_{j=1}^{k}g(b_j,c_j)}.
\end{eqnarray*}
As expected, $G(k)$ is indeed independent of $m$. Thus, (\ref{111})
holds. Inserting $G(k)$ into (\ref{111}) and making some
simplifications, we get the desired result at once.

{\bf Remark.} Note that the constant $C$ in Lemma \ref{13yl} is not
unique. If we define
$$
\prod^{0}_{j=1}f(a_j,c_j)g(b_j,d_j) =1,
\prod^{-1}_{j=1}f(a_j,c_j)g(b_j,d_j)=0,
$$
then we find that $C=-1$. If we redefine
$$
\prod^{0}_{j=1}f(a_j,c_j)g(b_j,d_j) =1,
\prod^{-1}_{j=1}f(a_j,c_j)g(b_j,d_j)=\frac{1}{f(a_0,c_0)g(b_0,d_0)},
$$
then
$$C=\frac{f(a_0,d_0)g(c_0,b_0)}{f(a_0,c_0)g(b_0,d_0)}.$$
To avoid this confusion, also with an effort to extend this result
to the setting of bilateral summation, we employ the convention of
defining (cf. \cite[Eq.(3.6.12)]{rahman})
\begin{eqnarray*}
\prod_{j=k}^mA_j= \left\{\begin{array}{ll} A_{k}A_{k+1}\cdots
A_{m},& m\geq k;\\
1,& m=k-1;\\
 (A_{m+1}A_{m+2}\cdots A_{k-1})^{-1},& m\leq k-2
\end{array}
\right.
\end{eqnarray*}
over $Z$.

Now, we are in a position to show our main theorem which gives the
bilateral form of Lemma \ref{13yl}.
\begin{dl}\label{13dl} Let $\{a_i\},\{b_i\},\{c_i\},\{d_i\}$ be arbitrary
sequences such that none of the denominators in (\ref{news111555})
vanish. Then for any nonnegative integers $m,n$,
\begin{align}
\sum^m_{k=-n} f(a_k,b_k)g(c_k,d_k)\frac{\prod_{j=1}^{k-1}f(a_j,c_j)}
{\prod_{j=1}^{k} f(a_j,d_j)}\frac{\prod_{j=1}^{k-1}g(b_j,d_j)}
{\prod_{j=1}^{k}g(b_j,c_j)}\nonumber\\=\frac{\prod_{j=1}^{m}f(a_j,c_j)}
{\prod_{j=1}^{m}f(a_j,d_j)}\frac{\prod_{j=1}^{m} g(b_j,d_j)}
{\prod_{j=1}^{m}g(b_j,c_j)}-\frac{\prod^{0}_{j=-n}f(a_j,d_j)}
{\prod^{0}_{j=-n}f(a_j,c_j)}\frac{\prod^{0}_{j=-n} g(b_j,c_j)}
{\prod^{0}_{j=-n}g(b_j,d_j)}\label{news111555}
\end{align}
provided $f\bot g$.
\end{dl}
{\sl Proof.} Assume that
$$
\sum_{k=-n}^mt_k=\frac{\prod_{j=1}^{m}f(a_j,c_j)}
{\prod_{j=1}^{m}f(a_j,d_j)}\frac{\prod_{j=1}^{m} g(b_j,d_j)}
{\prod_{j=1}^{m}g(b_j,c_j)}-\frac{\prod^{0}_{j=-n}f(a_j,d_j)}
{\prod^{0}_{j=-n}f(a_j,c_j)}\frac{\prod^{0}_{j=-n} g(b_j,c_j)}
{\prod^{0}_{j=-n}g(b_j,d_j)}.
$$
Apply the same argument as in Lemma \ref{13yl} only to find
\begin{eqnarray*}t_k=
    f(a_k,b_k)g(c_k,d_k)\frac{\prod_{j=1}^{k-1}f(a_j,c_j)}
{\prod_{j=1}^{k} f(a_j,d_j)}\frac{\prod_{j=1}^{k-1}g(b_j,d_j)}
{\prod_{j=1}^{k}g(b_j,c_j)}
\end{eqnarray*}
for $k\geq 0$. However for $k\leq -1$,
\begin{eqnarray*}t_k=
       f(a_k,b_k)g(c_k,d_k)\frac{\prod_{j=k+1}^{0}f(a_j,d_j)}
{\prod_{j=k}^{0} f(a_j,c_j)}\frac{\prod_{j=k+1}^{0}g(b_j,c_j)}
{\prod_{j=k}^{0}g(b_j,d_j)}.
\end{eqnarray*}
In this instance, by the above convention,  $t_k$ can be still
rewrite as
\begin{eqnarray*}t_k=
    f(a_k,b_k)g(c_k,d_k)\frac{\prod_{j=1}^{k-1}f(a_j,c_j)}
{\prod_{j=1}^{k} f(a_j,d_j)}\frac{\prod_{j=1}^{k-1}g(b_j,d_j)}
{\prod_{j=1}^{k}g(b_j,c_j)}.
 \end{eqnarray*}
This yields
\begin{align*}
\sum^m_{k=-n} f(a_k,b_k)g(c_k,d_k)\frac{\prod_{j=1}^{k-1}f(a_j,c_j)}
{\prod_{j=1}^{k} f(a_j,d_j)}\frac{\prod_{j=1}^{k-1}g(b_j,d_j)}
{\prod_{j=1}^{k}g(b_j,c_j)}\nonumber\\=\frac{\prod_{j=1}^{m}f(a_j,c_j)}
{\prod_{j=1}^{m}f(a_j,d_j)}\frac{\prod_{j=1}^{m} g(b_j,d_j)}
{\prod_{j=1}^{m}g(b_j,c_j)}-\frac{\prod^{0}_{j=-n}f(a_j,d_j)}
{\prod^{0}_{j=-n}f(a_j,c_j)}\frac{\prod^{0}_{j=-n} g(b_j,c_j)}
{\prod^{0}_{j=-n}g(b_j,d_j)}+C,
\end{align*}
where $C$ is independent  of $m,n$. When $m=n=0$, then from $f\bot
g$ it follows $C=0$. This gives the complete proof of
(\ref{news111555}).

 From now on, we call  (\ref{news111555}) the $(f,g)$-summation formula.
As we will see later, the $(f,g)$-summation formula indeed unifies
and extends all known  bibasic  summation formulae. It should be
pointed that the corresponding proofs  of these summation formulas
seem somewhat mysterious because they depend heavily on a tricky
factorization
 of a difference of two four-term
products into a four-term product (cf.\cite[Eq.(3.6.10)]{rahman}).
In the author's view, the $(f,g)$-summation formula provides a
natural interpretation about such ``mysterious" phenomena. Before we
turn to illustrate this, it is worth noting two particular cases of
the $(f,g)$-summation formula.

Set $c_j=b_0, d_j=x,n=0$ in (\ref{news111555}). Then the
$(f,g)$-summation formula reduces to
\begin{tl}\label{1.3dl} With the assumption as in Theorem \ref{13dl}. Then for any integer $m\geq 0$,
\begin{align}
\sum^m_{k=0}\frac{f(a_k,b_k)}{f(a_0,b_0)}\frac{\prod_{j=0}^{k-1}
f(a_j,b_0)}
{\prod_{j=1}^{k}g(b_j,b_0)}\frac{\prod_{j=0}^{k-1}g(b_j,x)}
{\prod_{j=1}^{k} f(a_j,x)}=\frac{\prod_{j=1}^{m}f(a_j,b_0)}
{\prod_{j=1}^{m}g(b_j,b_0)}\frac{\prod_{j=1}^{m} g(b_j,x)}
{\prod_{j=1}^{m}f(a_j,x)}\label{news11}
\end{align}
provided $f\bot g$. \end{tl}

Another immediate consequence of Theorem  \ref{13dl}  is Chu's
extension (cf.\cite[Theorem A]{chu}) of Gasper's bibasic summation
formula, which turned out to be useful in deriving
 hypergeometric identities. To state
Chu's result,  define
\begin{eqnarray*}
&&\phi^{\ast}(x;m)=\prod_{i=m}^{\infty}(a_i+xb_i); \psi^{\ast}(y;m)=\prod_{i=m}^{\infty}(c_i+yd_i);\\
&&\phi(x;m)=\prod_{i=0}^{m-1}(a_i+xb_i)=\phi^{\ast}(x;0)/\phi^{\ast}(x;m);\\
&&\psi(y;m)=\prod_{j=0}^{m-1}(c_j+yd_j)=\psi^{\ast}(y;0)/\psi^{\ast}(y;m).
\end{eqnarray*}
It is easy to check that the definitions of $\phi$ and $\psi$ are
confirmed with the previous convention. The next corollary is just
Chu's result.
\begin{tl} Let $\phi$ and $\psi$ be given as above,
$m,n\in Z$.
\begin{eqnarray}
&(x-y)\sum_{k=m}^n(a_kd_k-b_kc_k)
\frac{\phi(x;k)\psi(y;k)}{\phi(y;k+1)\psi(x;k+1)}\nonumber\\
&=\frac{\phi(x;m)\psi(y;m)}{\phi(y;m)\psi(x;m)}-\frac{\phi(x;n+1)\psi(y;n+1)}{\phi(y;n+1)\psi(x;n+1)}.
\end{eqnarray}
\end{tl}
{\sl Proof.} Without any loss of generality, assume $m,n\geq 0$.
Take in (\ref{news1155}) $f(x,y)=g(x,y)=x-y$ and make the
substitution
\begin{eqnarray*}
&&a_{i}\rightarrow \frac{a_{i-1}}{b_{i-1}},b_{i}\rightarrow \frac{c_{i-1}}{d_{i-1}};\\
&&c_{i}\rightarrow -x,d_{i}\rightarrow -y,\quad i=0,1,2,\cdots
\end{eqnarray*}
So the products are
$$
\prod_{j=1}^{k} f(a_j,c_j)=\frac{\phi(x;k)}{b_0b_1\cdots
b_{k-1}},\,\,\prod_{j=1}^{k}
g(b_j,d_j)=\frac{\psi(y;k)}{d_0d_1\cdots d_{k-1}},
$$
and
\begin{eqnarray*}
 \frac{\prod_{j=1}^{k-1}f(a_j,c_j)}
{\prod_{j=1}^{k}
f(a_j,d_j)}=\frac{b_{k-1}\phi(x;k-1)}{\phi(y;k)},\,\,\frac{\prod_{j=1}^{k-1}g(b_j,d_j)}
{\prod_{j=1}^{k}g(b_j,c_j)}=\frac{d_{k-1}\psi(y;k-1)}{\psi(x;k)}.
 \end{eqnarray*}
 Thus, (\ref{news1155}) becomes
\begin{eqnarray*}
&&(y-x)\sum_{k=0}^{n-1}(a_{k-1}d_{k-1}-b_{k-1}c_{k-1})
\frac{\phi(x;k-1)\psi(y;k-1)}{\phi(y;k)\psi(x;k)}=\frac{\phi(x;n)\psi(y;n)}{\phi(y;n)\psi(x;n)}+C.
\end{eqnarray*}
Substitute $n$ by $n+1$ and $k$ by $k+1$. This identity can be
reformulated as
\begin{eqnarray}
&&(y-x)\sum_{k=0}^{n}(a_{k}d_{k}-b_{k}c_{k})
\frac{\phi(x;k)\psi(y;k)}{\phi(y;k+1)\psi(x;k+1)}=\frac{\phi(x;n+1)\psi(y;n+1)}{\phi(y;n+1)\psi(x;n+1)}+C.\nonumber\\
&&\label{rr}
\end{eqnarray}
Replace $n$ by $m-1$. So (\ref{rr}) becomes
\begin{eqnarray}
&&(y-x)\sum_{k=0}^{m-1}(a_{k}d_{k}-b_{k}c_{k})
\frac{\phi(x;k)\psi(y;k)}{\phi(y;k+1)\psi(x;k+1)}=\frac{\phi(x;m)\psi(y;m)}{\phi(y;m)\psi(x;m)}+C.\nonumber\\
&&\label{ss}
\end{eqnarray}
Subtracting (\ref{ss}) from (\ref{rr}) leads to the desired result.

 Next,
we will exhibit some remarkable bibasic summation formulas from the
$(f,g)$-summation formula (\ref{news111555}) by specializing $f$ and
$g$ which are orthogonal to each other and related parameters. Among
them are Subarao and Verma's summation formula \cite{new1}, Gasper
and Rahman's bibasic summation formula \cite{8}, Chu's formula
\cite{chu} and Macdonald' extensions \cite{macd}. These facts show
convincingly that $\mbox{\sl Ker}\mathcal{L}^{(g)}_{3}$ provides
indeed a rich source of bibasic summation formulae.

{\bf I. $\left(1,x-y\right)$-summation formula.}
\begin{tl} With the assumption as in Theorem \ref{13dl}. Then for any integers $m,n\geq 0$,
\begin{align}
&&\sum^m_{k=-n} (c_k-d_k)\frac{\prod_{j=1}^{k-1}(b_j-d_j)}
{\prod_{j=1}^{k}(b_j-c_j)}=\frac{\prod_{j=1}^{m} (b_j-d_j)}
{\prod_{j=1}^{m}(b_j-c_j)}-\frac{\prod^{0}_{j=-n} (b_j-c_j)}
{\prod^{0}_{j=-n}(b_j-d_j)}.\label{news1133}
\end{align}
\end{tl}
This identity follows from Theorem \ref{13dl} by setting
\begin{eqnarray*}
  f(x,y)=1,g(x,y)=x-y,
\end{eqnarray*}
and the simple fact that $f\bot g$. It contains the following
Subarao and Verma's summation formula (cf. \cite[Eq.(3.1)]{new1}) as
a special case.
\begin{lz}
\begin{align}
\sum^m_{k=-n}az_k
\left(1-\frac{x_k}{az_k}\right)\left(1-\frac{y_k}{az_k}\right)
\frac{\prod_{j=1}^{k-1}(1-x_j)(1-y_j)}{\prod_{j=1}^k(1-a_jz_j)(1-\frac{x_jy_j}{az_j})}\nonumber\\
=\frac{\prod_{j=1}^{m}(1-x_j)(1-y_j)}{\prod_{j=1}^m(1-az_j)(1-\frac{x_jy_j}{az_j})}-
\frac{\prod_{j=-n}^0(1-az_j)(1-\frac{x_jy_j}{az_j})}{\prod_{j==-n}^{0}(1-x_j)(1-y_j)}.\label{news22}
\end{align}
\end{lz}
{\sl Proof.} Actually, make the substitution at first in
(\ref{news1133}) $b_j\mapsto t_j,c_j\mapsto
\frac{t_j-1}{1-y_j},d_j\mapsto\frac{t_j-1}{1-u_j}$ , and then
specificize parameters
$
  u_j\mapsto az_j,t_j\mapsto x_j/u_j.
$ It is easy to check by direct calculation that
\begin{align*}
&&\sum^m_{k=0} (c_k-d_k)\frac{\prod_{j=1}^{k-1}(b_j-d_j)}
 {\prod_{j=1}^{k}
(b_j-c_j)}=\sum^m_{k=-n}
\frac{(1-t_k)(u_k-y_k)}{(1-u_k)(1-y_k)}\frac{1-y_k}{1-t_ky_k}\prod_{j=1}^{k-1}\frac{(b_j-d_j)}
 {(b_j-c_j)}\\
 &&=\sum^m_{k=-n}
\frac{(1-t_k)(u_k-y_k)}{(1-u_k)(1-t_ky_k)}\prod_{j=1}^{k-1}
\frac{(1-x_j)(1-y_j)}{(1-az_j)(1-\frac{x_jy_j}{az_j})}.
\end{align*}
Note that
\begin{align*}
(1-t_k)(u_k-y_k)=az_k\left(1-\frac{x_k}{az_k}\right)\left(1-\frac{y_k}{az_k}\right);\\
(1-u_k)(1-t_ky_k)=(1-az_k)(1-\frac{x_ky_k}{az_k}).
\end{align*}
Thus, the right-hand side of the previous identity can be simplified
to
\begin{align*}
\sum^m_{k=-n}
az_k\left(1-\frac{x_k}{az_k}\right)\left(1-\frac{y_k}{az_k}\right)\,
\frac{\prod_{j=1}^{k-1}(1-x_j)(1-y_j)}{\prod_{j=1}^k(1-b_jz_j)(1-\frac{x_jy_j}{az_j})}
.
\end{align*}
This gives the desired result.

 {\bf II. $\left(x-y,x-y\right)$-summation
formula.}
\begin{tl}\label{coroll6} With the assumption as in Theorem \ref{13dl}. Then for any integers $m,n\geq 0$,
\begin{align}
\sum^m_{k=-n} (a_k-b_k)(c_k-d_k)\frac{\prod_{j=1}^{k-1}(a_j-c_j)}
{\prod_{j=1}^{k} (a_j-d_j)}\frac{\prod_{j=1}^{k-1}(b_j-d_j)}
{\prod_{j=1}^{k}(b_j-c_j)}\nonumber\\
=\frac{\prod_{j=1}^{m}(a_j-c_j)}
{\prod_{j=1}^{m}(a_j-d_j)}\frac{\prod_{j=1}^{m} (b_j-d_j)}
{\prod_{j=1}^{m}(b_j-c_j)}-\frac
{\prod_{j=-n}^{0}(a_j-d_j)}{\prod_{j=-n}^{0}(a_j-c_j)}\frac
{\prod_{j=-n}^{0}(b_j-c_j)}{\prod_{j=-n}^{0}(b_j-d_j)}.\label{news1122}
\end{align}
\end{tl}
This identity follows from Theorem \ref{13dl} by setting
\begin{eqnarray*}
  f(x,y)=x-y,g(x,y)=x-y.
\end{eqnarray*}
By Corollary \ref{coll21}, we see that $f\bot g$. It is worth
mentioning that this result contains Subarao and Verma's summation
formula \cite[Eq.(2.1)]{new1} with four independent bases. We
restate it as follows:
\begin{lz} Let $\{u_i\},\{v_i\},\{w_i\},\{z_i\}$ be arbitrary
sequences such that none of the denominators in (\ref{1.4}) vanish,
$m,n$ be nonnegative integers. Then the following holds,
\begin{align}
\sum^m_{k=-n}\frac{u_kv_kw_k}{z_k}
\left(1-u_kv_kw_kz_k\right)\left(1-\frac{w_kz_k}{u_kv_k} \right)
\left(1-\frac{u_kz_k}{ v_kw_k} \right) \left(
1-\frac{v_kz_k}{u_kw_k}\right)\nonumber\\
\frac{\prod_{j=1}^{k-1}(1-u_j^2)(1-v_j^2)(1-w_j^2)(1-z_j^2)}{\prod_{j=1}^k(1-\frac{u_jv_jw_j}{z_j})(1-\frac{u_jv_jz_j}{w_j})(1-\frac{w_jz_ju_j}{v_j})
(1-\frac{w_jz_jv_j}{u_j})}
\nonumber\\
=\prod_{j=1}^m\frac{(1-u_j^2)(1-v_j^2)(1-w_j^2)(1-z_j^2)}{(1-\frac{u_jv_jw_j}{z_j})(1-\frac{u_jv_jz_j}{w_j})(1-\frac{w_jz_ju_j}{v_j})
(1-\frac{w_jz_jv_j}{u_j})}\nonumber\\
-\prod_{j=-n}^{0}
\frac{(1-\frac{u_jv_jw_j}{z_j})(1-\frac{u_jv_jz_j}{w_j})(1-\frac{w_jz_ju_j}{v_j})
(1-\frac{w_jz_jv_j}{u_j})}{(1-u_j^2)(1-v_j^2)(1-w_j^2)(1-z_j^2)}.\label{1.4}
\end{align}
\end{lz}
 {\sl Proof.}
Specify all parameters in (\ref{news1122}) by
\begin{eqnarray*}
  a_j &=& 1/A_{2,j},b_j=1/B_{2,j}; \\
  c_j &=& Y_j,d_j=X_j ,
 \end{eqnarray*}
 and then set parameters in the corresponding result by
 $$A_{2,j}Y_j=A_{1,j},B_{2,j}X_j=B_{1,j}.$$ Then, after some calculation, we arrive at
\begin{align}
\sum^m_{k=-n} (1/A_{2,k}-1/B_{2,k})
(\frac{A_{1,k}}{A_{2,k}}-\frac{B_{1,k}}{B_{2,k}})A_{2,k}B_{2,k}\nonumber\\
\frac {\prod_{j=1}^{k-1}(1-A_{1,j})}{\prod_{j=1}^{k}
(1-A_{2,j}\frac{B_{1,j}}{B_{2,j}})} \frac{\prod_{j=1}^{k-1}
(1-B_{1,j})}{\prod_{j=1}^{k}(1-B_{2,j}\frac{A_{1,j}}{A_{2,j}})}
=\frac{\prod_{j=1}^{m}(1-A_{1,j})}{\prod_{j=1}^{m}(1-A_{2,j}\frac{B_{1,j}}{B_{2,j}})}
\frac{\prod_{j=1}^{m}(1-B_{1,j})}{\prod_{j=1}^{m}
(1-B_{2,j}\frac{A_{1,j}}{A_{2,j}})}-C_n,\label{b22}
\end{align}
where the term $C_n$ is
$$
C_n=\prod_{j=-n}^{0}\frac{(1-A_{2,j}\frac{B_{1,j}}{B_{2,j}})(1-B_{2,j}\frac{A_{1,j}}{A_{2,j}})}
{(1-A_{1,j})(1-B_{1,j})}.
$$

 In order to show Subarao and Verma's
result (\ref{1.4}), we  make the further substitution of parameters
\begin{eqnarray*}
  A_{1,j} &=& \frac{u_j^2+v_j^2}{1+u_j^2v_j^2},
  A_{2,j}= \frac{u_jv_j}{1+u_j^2v_j^2}; \\
 B_{1,j} &=& \frac{w_j^2+z_j^2}{1+w_j^2z_j^2};
 B_{2,j}=\frac{w_jz_j}{1+w_j^2z_j^2},
\end{eqnarray*}
which in turn gives
\begin{eqnarray*}
  A_{1,j}/A_{2,j} &=& u_j/v_j+v_j/u_j \\
  B_{1,j}/B_{2,j} &=& w_j/z_j+z_j/w_j.
\end{eqnarray*}
Using these relations to calculate
\begin{eqnarray*}
  \frac {1-A_{1,j}}{
1-A_{2,j}\frac{B_{1,j}}{B_{2,j}}}
=\frac{(1-u_j^2)(1-v_j^2)}{(1-\frac{u_jv_jw_j}{z_j})(1-\frac{u_jv_jz_j}{w_j})}.
\end{eqnarray*}
Also, similar calculation gives that
\begin{eqnarray*}
&&  \frac {1-B_{1,j}}{1-B_{2,j}\frac{A_{1,j}}{A_{2,j}}}=
\frac{(1-w_j^2)(1-z_j^2)}{(1-\frac{w_jz_ju_j}{v_j})
(1-\frac{w_jz_jv_j}{u_j})};\\
&&\frac{A_{2,k}B_{2,k}}{(1-A_{2,k}\frac{B_{1,k}}
{B_{2,k}})(1-B_{2,k}\frac{A_{1,k}}{A_{2,k}})}=\frac{u_kv_kw_kz_k}{(1-\frac{u_kv_kw_k}{z_k})(1-\frac{u_kv_kz_k}{w_k})
(1-\frac{w_kz_ku_k}{v_k})(1-\frac{w_kz_kv_k}{u_k})}.\\
\end{eqnarray*}
So, the production is
\begin{eqnarray*}
&& A_{2,k}B_{2,k}\frac
{\prod_{j=1}^{k-1}(1-A_{1,j})}{\prod_{j=1}^{k}
(1-A_{2,j}\frac{B_{1,j}}{B_{2,j}})} \frac{\prod_{j=1}^{k-1}
(1-B_{1,j})}{\prod_{j=1}^{k}(1-B_{2,j}\frac{A_{1,j}}{A_{2,j}})}\\
&&=u_kv_kw_kz_k\frac
{\prod_{j=1}^{k-1}(1-u_j^2)(1-v_j^2)(1-w_j^2)(1-z_j^2)}{\prod_{j=1}^{k}(1-\frac{u_jv_jw_j}{z_j})(1-\frac{u_jv_jz_j}{w_j})
(1-\frac{w_jz_ju_j}{v_j})(1-\frac{w_jz_jv_j}{u_j})}.
\end{eqnarray*}
Further calculation yields
\begin{eqnarray*}
&&
\frac{1}{A_{2,k}}-\frac{1}{B_{2,k}}=\frac{1+u_k^2v_k^2}{u_kv_k}-\frac{1+w_k^2z_k^2}{w_kz_k}
= \frac{\left(u_kv_kw_kz_k-1 \right) \left( v_ku_k-z_kw_k \right)
}{u_kv_kw_kz_k};
\\
&&\frac{A_{1,k}}{A_{2,k}}-\frac{B_{1,k}}{B_{2,k}}
=\frac{u_k}{v_k}+\frac{v_k}{u_k}-\frac{w_k}{z_k}-\frac{z_k}{w_k}=
\frac{\left( v_kw_k-z_ku_k \right)  \left(-u_kw_k+v_kz_k \right)
}{u_kv_kw_kz_k}.
\end{eqnarray*}
This also gives
\begin{eqnarray*}
  (\frac{1}{A_{2,k}}-\frac{1}{B_{2,k}})(\frac{A_{1,k}}{A_{2,k}}-
  \frac{B_{1,k}}{B_{2,k}})\\
  =
  \frac{\left(1-u_kv_kw_kz_k\right)
\left( u_kv_k-w_kz_k \right) \left( v_kw_k-u_kz_k \right)  \left(
u_kw_k-v_kz_k \right) }{(u_kv_kw_kz_k)^2}.
 \end{eqnarray*}
 Insert all these into (\ref{b22}) and rearrange the resulting identity. The result is (\ref{1.4}).

Observe that by performing various substitutions we may yet deduce
other summation formulas obtained by Subarao and Verma \cite{new1}.
We leave them for the interested reader.

Further, making the substitution $n\mapsto 0, b_j\mapsto
b_0/b_j,c_j\mapsto 1,d_j\mapsto x$ in Corollary \ref{coroll6} gives
Chu's extension of Krattenthaler's matrix inversion \cite{chu}.

\begin{lz} Let $\{a_i\},\{b_i\}$ be arbitrary
sequences, $x$ be indeterminate. Then for any integer $m\geq 0$,
\begin{align*}
&&\sum^m_{k=0}
\frac{b_0-b_ka_k}{b_0-b_0a_0}\frac{\prod_{j=0}^{k-1}(1-a_j)(b_0-b_jx)}
{\prod_{j=1}^{k}(1-a_j/x)(b_0-b_j)}\frac{1}{x^k}=\frac{\prod_{j=1}^{m}(1-a_j)(b_0-b_jx)}
{\prod_{j=1}^{m}(1-a_j/x)(b_0-b_j)}\frac{1}{x^m}.
\end{align*}
\end{lz}
 {\bf III.\,$\left((1-axy)(1-b\frac{x}{y}),(x-y)(1-\frac{b}{axy})\right)$-summation formula.}\\
 Define that
\begin{eqnarray*}
  g(x,y)=(x-y)(1-\frac{b}{axy}),\,\, f(x,y)=(1-axy)(1-b\frac{x}{y}).
\end{eqnarray*}
Then by Corollary \ref{coll23}, $f\bot g$. This reduces Theorem
\ref{13dl} to
\begin{tl} With the assumption as in Theorem \ref{13dl}. Then for any integer $m,n\geq 0$,
\begin{align}
\sum^m_{k=-n}
(1-aa_kb_k)(1-b\frac{a_k}{b_k})(c_k-d_k)(1-\frac{b}{ac_kd_k})\nonumber\\
\times\frac{\prod_{j=1}^{k-1}(1-aa_jc_j)(1-b\frac{a_j}{c_j})}
{\prod_{j=1}^{k} (1-aa_jd_j)(1-b\frac{a_j}{d_j})}
\frac{\prod_{j=1}^{k-1}(b_j-d_j)(1-\frac{b}{ab_jd_j}) }
{\prod_{j=1}^{k}(b_j-c_j)(1-\frac{b}{ab_jc_j})}\nonumber\\
=\frac{\prod_{j=1}^{m}(1-aa_jc_j)(1-b\frac{a_j}{c_j})}
{\prod_{j=1}^{m} (1-aa_jd_j)(1-b\frac{a_j}{d_j})}
\frac{\prod_{j=1}^{m}(b_j-d_j)(1-\frac{b}{ab_jd_j}) }
{\prod_{j=1}^{m}(b_j-c_j)(1-\frac{b}{ab_jc_j})}\nonumber\\- \frac
{\prod_{j=-n}^{0}
(1-aa_jd_j)(1-b\frac{a_j}{d_j})}{\prod_{j=-n}^{0}(1-aa_jc_j)(1-b\frac{a_j}{c_j})}
\frac{\prod_{j=-n}^{0}(b_j-c_j)(1-\frac{b}{ab_jc_j})}{\prod_{j=-n}^{0}(b_j-d_j)(1-\frac{b}{ab_jd_j})
}\label{00}.
\end{align}
\end{tl}
Identity (\ref{00}) contains Gasper and Rahman's indefinite bibasic
summation formula \cite[Eq.(3.6.3)]{10} as a special case. In fact,
set
\begin{eqnarray*}
  a_k=p^k, b_k=dq^k, c_k=1, d_k=\frac{d}{x}.
\end{eqnarray*}
 It is easily seen that
\begin{align}
\sum_{k=-n}^{m}\frac{(1-adp^kq^{k})(1-b/dp^kq^{-k})} {(1-ad)(1-b/d)}
\frac{(a,b;p)_k(x,ad^2/(bx);q)_k}{(dq,adq/b;q)_k(adp/x,bpx/d;p)_k}q^k\nonumber\\
=\frac{(1-a)(1-b)(1-x)(1-ad^2/(bx))} {(1-ad)(1-b/d)(d-x)(1-ad/(bx))}
\left\{\frac{(ap,bp;p)_m(xq,ad^2q/(bx);q)_m}
{(dq,adq/b;q)_m(adp/x,bpx/d;p)_m}\right.\nonumber\\
-\left.\frac{(x/(ad),d/(bx);p)_{n+1}(1/d,b/(ad);q)_{n+1}}
{(1/x,bx/(ad^2);q)_{n+1}(1/a,1/b;p)_{n+1}}\right\}.\label{newf}
\end{align}
Gasper and Rahman used (\ref{newf}) to set up a series of quadratic
and cubic summation and transformation formulas of the basic
hypergeometric series. See \cite[Section 3.8]{rahman,10} for more
details.

Further, let $n=0,b=0,d=1$ or $n=0,d=1$. In both cases,
$(1/d;q)_{n+1}=0$. After some simplification, this identity reduces
to Gosper's bibasic summation formula
\begin{eqnarray}   & &\sum^m_{k=0} \frac{1-ap^kq^k}{1-a} \,
\frac{(a;p)_k(1/x;q)_k}{(q;q)_k(apx;p)_k}\, x^k
=\frac{(ap;p)_m(q/x;q)_m}{(q;q)_m(apx;p)_m}\, x^{m}\label{4.7}
\end{eqnarray}
 and  Gasper's bibasic summation formula
\begin{align}
\sum_{k=0}^{m}\frac{(1-ap^kq^{k})(1-bp^kq^{-k})}{(1-a)(1-b)}
\frac{(a,b;p)_k(x,a/(bx);q)_k}{(q,aq/b;q)_k(ap/x,bpx;p)_k}q^k\nonumber\\
=\frac{(ap,bp;p)_m(xq,aq/(bx);q)_m}{(q,aq/b;q)_m
(ap/x,bpx;p)_m},\label{newff}
\end{align}
respectively.

{\bf IV. $\left((x+y)(x+\frac{b}
  {ay}),(x-y)(1-\frac{b}{axy})\right)$-summation
formula.}\\  Define that
\begin{eqnarray*}
  g(x,y)=(x-y)(1-\frac{b}{axy}),\,\,
  f(x,y)=(x+y)(x+\frac{b}
  {ay}).
  \end{eqnarray*}
Then  $f\bot g$. Therefore, Theorem \ref{1.3dl} reduces to
\begin{tl} With the assumption as in Theorem \ref{13dl}. Then for any integers $m,n\geq 0$,
\begin{align}
\sum_{k=-n}^{m}
(a_k+b_k)(a_k+\frac{b}{ab_k})(c_k-d_k)(1-\frac{b}{ac_kd_k})\nonumber\\
\times\frac{\prod_{j=1}^{k-1}(a_j-c_j)(a_j+\frac{b}{ac_j})}
{\prod_{j=1}^{k} (a_j+d_j)(a_j+\frac{b}{ad_j})}
\frac{\prod_{j=1}^{k-1}(b_j-d_j)(1-\frac{b}{ab_jd_j}) }
{\prod_{j=1}^{k}(b_j-c_j)(1-\frac{b}{ab_jc_j})}\nonumber\\
=\frac{\prod_{j=1}^{m}(a_j-c_j)(a_j+\frac{b}{ac_j})}
{\prod_{j=1}^{m} (a_j+d_j)(a_j+\frac{b}{ad_j})}
\frac{\prod_{j=1}^{m}(b_j-d_j)(1-\frac{b}{ab_jd_j}) }
{\prod_{j=1}^{m}(b_j-c_j)(1-\frac{b}{ab_jc_j})}\nonumber\\
-\frac{\prod_{j=-n}^{0}
(a_j+d_j)(a_j+\frac{b}{ad_j})}{\prod_{j=-n}^{0}(a_j-c_j)(a_j+\frac{b}{ac_j})}
\frac{\prod_{j=-n}^{0}(b_j-c_j)(1-\frac{b}{ab_jc_j})}{\prod_{j=-n}^{0}(b_j-d_j)(1-\frac{b}{ab_jd_j})
} .
\end{align}
\end{tl}
{\bf V.
$\left(y(1-\frac{x}{y})(1-\frac{xy}{d}),y(1-\frac{x}{y})(1-\frac{xy}{d})\right)$-summation
formula.}

 Define that
\begin{eqnarray*}
  g(x,y)=f(x,y)=y(1-x/y)(1-xy/d).
  \end{eqnarray*}
Corollary \ref{coll21} states that $f\bot g$. In this case, Theorem
\ref{1.3dl} reduces to
\begin{tl} \label{gaoxin}With the assumption as in Theorem \ref{13dl}.  Then the following holds,
\begin{align}
\sum_{k=-n}^{m}
(b_k-a_k)(1-\frac{a_kb_k}{d})(d_k-c_k)(1-\frac{c_kd_k}{d})\nonumber\\
\times\frac{\prod_{j=1}^{k-1}(c_j-a_j)(1-\frac{a_jc_j}{d})}
{\prod_{j=1}^{k} (d_j-a_j)(1-\frac{a_jd_j}{d})}
\frac{\prod_{j=1}^{k-1}(d_j-b_j)(1-\frac{b_jd_j}{d}) }
{\prod_{j=1}^{k}(c_j-b_j)(1-\frac{b_jc_j}{d})}\nonumber\\
=\frac{\prod_{j=1}^{m}(c_j-a_j)(1-\frac{a_jc_j}{d})}
{\prod_{j=1}^{m} (d_j-a_j)(1-\frac{a_jd_j}{d})}
\frac{\prod_{j=1}^{m}(d_j-b_j)(1-\frac{b_jd_j}{d}) }
{\prod_{j=1}^{m}(c_j-b_j)(1-\frac{b_jc_j}{d})}\nonumber\\
-\frac{\prod_{j=-n}^{0}
(d_j-a_j)(1-\frac{a_jd_j}{d})}{\prod_{j=-n}^{0}(c_j-a_j)(1-\frac{a_jc_j}{d})}
\frac{\prod_{j=-n}^{0}(c_j-b_j)(1-\frac{b_jc_j}{d})}{\prod_{j=-n}^{0}(d_j-b_j)(1-\frac{b_jd_j}{d})
} .\label{news}
\end{align}
\end{tl}

An important case of this identity is Chu's generalization of Gasper
and Rahman's formula, i.e., (\ref{newf}).
\begin{lz}\label{chu} {\rm (cf.\cite{chu})} With the assumption as Corollary \ref{gaoxin}. The following holds
\begin{align}
\sum_{k=0}^{m}(b_0-a_kb_k)(b_k-a_k
b_0/d)\frac{\prod_{j=1}^{k-1}(1-a_j)(1-a_j/d)(b_0- b_jx)(b_0 - b_j
d/x)} {\prod_{j=1}^{k}
(b_0-b_j)(b_0-b_jd)(1-a_j/x)(1-a_jx/d)}\nonumber\\
=\frac{x}{(d-x)(x-1)}\prod_{j=1}^{m}\frac{(1-a_j)(1-a_j/d)(b_0-
b_jx)(b_0 - b_j d/x)} {(b_0-b_j)(b_0-b_jd)(1-a_j/x)(1-a_j x/d)}.
\end{align}
\end{lz}
{\sl Proof.} Set in (\ref{news}) $n\mapsto 0, b_i\mapsto
b_id/b_0=b^{'}_i,c_i\mapsto 1,d_i\mapsto x$. It is easy to check
that $b_0^{'}=d$, and
$$f(a_k,b^{'}_k)= (b_kd/b_0)(b_0-a_kb_k)(b_k-a_k
b_0/d)/(b_0b_k)=d(b_0-a_kb_k)(b_k-a_k b_0/d)/b_0^2,$$ as well as
$$\frac{\prod_{j=1}^{k-1}g(b^{'}_j,x)}{\prod_{j=1}^{k}f(a_j,x)}=\frac{1}{xb_0^{2k-2}}\frac{\prod_{j=1}^{k-1}
(b_0- b_jx)(b_0 - b_j d/x)} {\prod_{j=1}^{k} (1-a_j/x)(1-a_jx/d)}.$$
To simplify (\ref{news}) by these calculation gives the desired
result.

Apply  the substitution  $b_i\mapsto deb_i,c_i\mapsto 1$, and
$d_i\mapsto x/e$ to (\ref{news}). Then we get
\begin{lz}\label{chu} {\rm (cf.\cite[Eq.(4.32)]{macd})} With the assumption as Corollary \ref{gaoxin}. The following holds
\begin{align}
\sum_{k=0}^{m}
(b_k-a_k/(de))(1-ea_kb_k)\nonumber\\
\times\frac{\prod_{j=1}^{k-1}(1-a_j)(1-a_j/d)} {\prod_{j=1}^{k}
(1-a_je/x)(1-\frac{a_jx}{de})}
\frac{\prod_{j=1}^{k-1}(1-b_jde^2/x)(1-xb_j) }
{\prod_{j=1}^{k}(1-deb_j)(1-eb_j)}\nonumber\\
=\frac{x}{(x-e)(x-de)}
\left\{\frac{(1-a_0e/x)(1-\frac{a_0x}{de})(1-deb_0)(1-eb_0)}
{(1-b_0de^2/x)(1-xb_0)(1-a_0)(1-a_0/d)}\right.\nonumber\\
 \left.-\frac{\prod_{j=1}^{m}(1-a_j)(1-a_j/d)} {\prod_{j=1}^{m}
(1-a_je/x)(1-\frac{a_jx}{de})}
\frac{\prod_{j=1}^{m}(1-b_jde^2/x)(1-xb_j)}
{\prod_{j=1}^{m}(1-deb_j)(1-eb_j)}\right\}.
\end{align}
\end{lz}

In his private communication with G.Bhatnagar, Macdonald extended
this result further to the following  \cite[Theorems 2.27 and 2.29,
p.200-201]{macd}
\begin{lz}
 Let $\{a_i\},\{b_i\},\{c_i\},\{d_i\}$ be arbitrary sequences
, $e$ be an indeterminate, $m,n$ be nonnegative integers, such that
none of the denominators in (\ref{mac}) vanish. Then
\begin{align}
\sum_{k=-n}^{m} e(1-a_kb_ke)(b_k-a_k/(d_ke))(1-c_k/e)(1-d_ke/c_k)
\nonumber\\
\times\frac{\prod_{j=1}^{k-1}(1-a_j) (1-\frac{a_j}{d_j})}
{\prod_{j=1}^{k} (1-\frac{a_je}{c_j})(1-\frac{a_jc}{d_je})}
\frac{\prod_{j=1}^{k-1}(1-b_jc)(1-\frac{b_jd_je^2}{c_j})}
{\prod_{j=1}^{k} (1-b_je)(1-b_jd_je)}\nonumber\\
=\frac{\prod_{j=1}^{m}(1-a_j) (1-\frac{a_j}{d_j})} {\prod_{j=1}^{m}
(1-\frac{a_je}{c_j})(1-\frac{a_jc_j}{d_je})}
\frac{\prod_{j=1}^{m}(1-b_jc_j)(1-\frac{b_jd_je^2}{c_j})}
{\prod_{j=1}^{m}(1-b_je)(1-b_jd_je)}\nonumber\\
-\frac {\prod_{j=-n}^{0}
(1-\frac{a_je}{c_j})(1-\frac{a_jc_j}{d_je})}{\prod_{j=-n}^{0}(1-a_j)
(1-\frac{a_j}{d_j})} \frac
{\prod_{j=-n}^{0}(1-b_je)(1-b_jd_je)}{\prod_{j=-n}^{0}(1-b_jc_j)(1-\frac{b_jd_je^2}{c_j})}\label{mac}.
\end{align}
\end{lz}
{\sl Proof.}\,\, As a key fact in Macdonald's proof, Eq.(2.24) in
\cite {macd} can be  replaced by  $f\bot g$. In fact,
 for arbitrary sequences $\{a_i\},\{b_i\},\{c_i\},\{d_i\}$,
\begin{eqnarray*}
f(a_k,c_k)g(b_k,d_k)- g(b_k,c_k)f(a_k,d_k)=f(a_k,b_k)g(c_k,d_k),
\end{eqnarray*}
where $f,g$ are given as Corollary \ref{gaoxin}. This identity
remains valid under the substitution  $b_i\mapsto deb_i,c_i\mapsto
1$ , and $d_i\mapsto c/e$. The result is
\begin{eqnarray*}
 f(a_k,1)g(deb_k,c/e)- g(deb_k,1)f(a_k,c/e)=f(a_k,deb_k)g(1,c/e).
\end{eqnarray*}
Since $c,d$ in this identity  are arbitrary (if so, $e$ is not free
since each denominators in (\ref{mac}) should not equal zero), we
set $c\mapsto c_k,d\mapsto d_k$. As there should be no confusion, we
still write $f(x,y)$ and $g(x,y)$ for the resulting functions. Thus
\begin{eqnarray*}
f(a_k,1)g(d_keb_k,c_k/e)-
g(d_keb_k,1)f(a_k,c_k/e)=f(a_k,d_keb_k)g(1,c_k/e).
\end{eqnarray*}
It can be reformulated, under the assumption that
$g(d_keb_k,1)f(a_k,c_k/e)\neq 0$,  as
\begin{eqnarray*}
\frac{f(a_k,1)g(d_keb_k,c_k/e)}
{g(d_keb_k,1)f(a_k,c_k/e)}-1=\frac{f(a_k,d_keb_k)g(1,c_k/e)}
{g(d_keb_k,1)f(a_k,c_k/e)}.
\end{eqnarray*}
Based on this identity and proceed as in Theorem \ref{13dl}. Assume
again that
\begin{align*}
&&\sum_{k=-n}^{m}
w_k\frac{\prod_{j=1}^{k-1}f(a_j,d_jeb_j)g(1,c_j/e)}
{\prod_{j=1}^{k}g(d_jeb_j,1)f(a_j,c_j/e)}=\prod_{j=1}^{m}
\frac{f(a_j,1)g(d_jeb_j,c_j/e)} {g(d_jeb_j,1)f(a_j,c_j/e)}-C_n,
\end{align*}
where
$$C_n=\frac{\prod^{0}_{j=-n}f(a_j,d_j)}
{\prod^{0}_{j=-n}f(a_j,c_j)}\frac{\prod^{0}_{j=-n} g(b_j,c_j)}
{\prod^{0}_{j=-n}g(b_j,d_j)}.
$$
 Solving this identity to get
\begin{align*}
w_k=f(a_k,d_keb_k)g(1,c_k/e)
\frac{\prod_{j=1}^{k-1}f(a_j,1)g(d_jeb_j,c_j/e)}
{\prod_{j=1}^{k-1}f(a_j,d_jeb_j)g(1,c_j/e)}.
\end{align*}
It in turn reduces the previous identity to
\begin{align*}
\sum_{k=-n}^{m}
f(a_k,d_keb_k)g(1,c_k/e)\frac{\prod_{j=1}^{k-1}f(a_j,1)g(d_jeb_j,c_j/e)}
{\prod_{j=1}^{k}g(d_jeb_j,1)f(a_j,c_j/e)}\\=\prod_{j=1}^{m}
\frac{f(a_j,1)g(d_jeb_j,c_j/e)}
{g(d_jeb_j,1)f(a_j,c_j/e)}-\prod_{j=-n}^{0} \frac
{g(d_jeb_j,1)f(a_j,c_j/e)}{f(a_j,1)g(d_jeb_j,c_j/e)},
\end{align*}
 which can be simplified to the desired result.

{\bf VI.  $\left(f,f\right)$- elliptic summation formula.}\\
As an immediate consequence of (\ref{news11}) being combined with
the \emph{addition formula }of the elliptic function
(cf.\cite{1000}), i.e., Corollary \ref{coll22}, we  obtain the
following  identity of the elliptic hypergeometric series, which is
equivalent to after relabeling, as pointed out by the referee, (3.2)
of \cite{1000} and (3.6) of \cite{101}.
\begin{tl}
Let $\theta(x)$ be given as in Corollary \ref{coll22},
$\{a_i\},\{b_i\},\{c_i\},\{d_i\}$ be arbitrary sequences such that
none of the denominators in (\ref{news13}) vanish, $m,n$ be
nonnegative integers. Then the following holds,
\begin{align}
\sum_{k=-n}^{m}\frac{b_k}{c_k}\theta(a_kb_k)\theta(c_kd_k)\theta
\left(\frac{a_k}{b_k}\right)\theta\left(\frac{c_k}{d_k}\right)\frac{\prod_{j=1}^{k-1}
\theta(a_jc_j)\theta(a_j/c_j)}
{\prod_{j=1}^{k}\theta(b_jc_j)\theta(b_j/c_j)}\frac{\prod_{j=1}^{k-1}\theta(b_jd_j)\theta(b_j/d_j)}
{\prod_{j=1}^{k} \theta(a_jd_j)\theta(a_j/d_j)}\nonumber\\
=\frac{\prod_{j=1}^{m} \theta(a_jc_j)\theta(a_j/c_j)}
{\prod_{j=1}^{m}\theta(b_jc_j)\theta(b_j/c_j)}\frac{\prod_{j=1}^{m}\theta(b_jd_j)\theta(b_j/d_j)}
{\prod_{j=1}^{m} \theta(a_jd_j)\theta(a_j/d_j)}\nonumber\\
-\frac{\prod_{j=-n}^{0}\theta(b_jc_j)\theta(b_j/c_j)}{\prod_{j=-n}^{0}
\theta(a_jc_j)\theta(a_j/c_j)}\frac{\prod_{j=-n}^{0}
\theta(a_jd_j)\theta(a_j/d_j)}{\prod_{j=-n}^{0}\theta(b_jd_j)\theta(b_j/d_j)}\label{news13}.
\end{align}
\end{tl}
{\sl Proof.} By Corollary \ref{coll22}, we see that $f\bot f$. Then
the desired result is an immediate consequence of Theorem
\ref{13dl}.
\section{The $(f,g)$-inversions from the $(f,g)$-summation formula}
In \cite{0020}, the author introduced an operator in order to set up
the $(f,g)$-inversion. Some possibly approaches to show it in a
simple style were suggested by the anonymous referee(s) and tried by
the author.  So far, we have found it can be proved by
Krattenthaler's operator method \cite{18} (only  for the case $f=g$
) and induction  \cite{1000}. Surprisingly,  the $(f,g)$-summation
formula
 leads to an alternative proof of this matrix inversion, which is considerably
 simpler than one given in \cite{0020}. It is worth noting that the
similar argument appeared in Bhatnagar's proof of Krattenthaler's
matrix inversion \cite{macd1}.

 {\sl The proof of Theorem \ref{math4}.}\,\,
Set $c_j=b_0, d_j=x$ in (\ref{news111555}). Then the
$(f,g)$-summation formula reduces to
\begin{align}
\sum^n_{k=0}\frac{f(a_k,b_k)}{f(a_0,b_0)}\frac{\prod_{j=0}^{k-1}
f(a_j,b_0)}
{\prod_{j=1}^{k}g(b_j,b_0)}\frac{\prod_{j=0}^{k-1}g(b_j,x)}
{\prod_{j=1}^{k} f(a_j,x)}=\frac{\prod_{j=1}^{n}f(a_j,b_0)}
{\prod_{j=1}^{n}g(b_j,b_0)}\frac{\prod_{j=1}^{n} g(b_j,x)}
{\prod_{j=1}^{n}f(a_j,x)}.
\end{align}
Actually, this identity can also be reformulated as
\begin{align}
\sum^n_{k=0}f(a_k,b_k)\frac {\prod_{j=k+1}^{n}
f(a_j,x)}{\prod_{j=k}^{n}g(b_j,x)}\left\{g(b_0,x)\frac{\prod_{j=0}^{k-1}
f(a_j,b_0)}
{\prod_{j=1}^{k}g(b_j,b_0)}\right\}=\frac{\prod_{j=1}^{n}f(a_j,b_0)}
{\prod_{j=1}^{n}g(b_j,b_0)}.
\end{align}
Replace $n$ by $n-1$ and  set $x\mapsto b_n$. Then
\begin{align*}
\sum^{n-1}_{k=0}f(a_k,b_k)\frac {\prod_{j=k+1}^{n-1}
f(a_j,b_n)}{\prod_{j=k}^{n-1}g(b_j,b_n)}\left\{\frac{\prod_{j=0}^{k-1}
f(a_j,b_0)}
{\prod_{j=1}^{k}g(b_j,b_0)}\right\}=-\frac{\prod_{j=1}^{n-1}f(a_j,b_0)}
{\prod_{j=1}^{n}g(b_j,b_0)},
\end{align*}
which is
\begin{align*}
\sum^{n}_{k=0}f(a_k,b_k)\frac {\prod_{j=k+1}^{n-1}
f(a_j,b_n)}{\prod_{j=k}^{n-1}g(b_j,b_n)}\left\{\frac{\prod_{j=0}^{k-1}
f(a_j,b_0)} {\prod_{j=1}^{k}g(b_j,b_0)}\right\}=\delta_{0,n}
\end{align*}
On the other hand, it is easily seen that
\begin{align*}
\sum^{n}_{k=0}\frac{\prod_{j=k}^{n-1}f(a_j,b_k)}
{\prod_{j=k+1}^{n}g(b_j,b_k)}\delta_{0,k}=\left\{\frac{\prod_{j=0}^{n-1}
f(a_j,b_0)} {\prod_{j=1}^{n}g(b_j,b_0)}\right\}.
\end{align*}
 The result is that two matrices $F=(f_{n,k})$ and $G=(g_{n,k})$, whose
 entries given by
\begin{eqnarray*}
f_{n,k}&=&\frac{\prod_{i=k}^{n-1}f(a_i,b_k)}
{\prod_{i=k+1}^{n}g(b_i,b_k)}\label{a12}\qquad\mbox{and}\\
g_{n,k}&=&
\frac{f(a_k,b_k)}{f(a_n,b_n)}\frac{\prod_{i=k+1}^{n}f(a_i,b_n)}
{\prod_{i=k}^{n-1}g(b_i,b_n)},\qquad\mbox{respectively,}\label{b12}
\end{eqnarray*}
are inverses of each other. As desired.

In a similar technique, we can obtain the following identity
\begin{dl} Preserve the convention as before. Then for $m\geq 0, n\geq 1$,
\begin{align*}
\sum^{m}_{k=-n} f(a_k,b_k)
 \frac{\prod_{j=m}^{k-1}g(b_j,b_{m})}
{\prod_{j=m}^{k} f(a_j,b_{m})}\frac{\prod_{j=1}^{k-1}f(a_j,b_{-n})}
{\prod_{j=1}^{k}g(b_j,b_{-n})}=0.
\end{align*}
\end{dl}
{\sl Proof.} In (\ref{news111555}), set $c_j\mapsto b_{-n},
d_j\mapsto x$. Since $g(b_{-n},b_{-n})=0$, it gives
\begin{align}
\sum^m_{k=-n} f(a_k,b_k)g(
b_{-n},x)\frac{\prod_{j=1}^{k-1}f(a_j,b_{-n})} {\prod_{j=1}^{k}
f(a_j,x)}\frac{\prod_{j=1}^{k-1}g(b_j,x)}
{\prod_{j=1}^{k}g(b_j,b_{-n})}\nonumber\\=\frac{\prod_{j=1}^{m}f(a_j,b_{-n})}
{\prod_{j=1}^{m}f(a_j,x)}\frac{\prod_{j=1}^{m} g(b_j,x)}
{\prod_{j=1}^{m}g(b_j,b_{-n})}.
\end{align}
Set $x \mapsto b_{m+1}$ again in this identity. The result is
\begin{align*}
\sum^m_{k=-n} f(a_k,b_k)g(
b_{-n},b_{m+1})\frac{\prod_{j=1}^{k-1}f(a_j,b_{-n})}
{\prod_{j=1}^{k}g(b_j,b_{-n})}
\frac{\prod_{j=1}^{k-1}g(b_j,b_{m+1})}{\prod_{j=1}^{k}
f(a_j,b_{m+1})}
\\=\frac{\prod_{j=1}^{m}f(a_j,b_{-n})}{\prod_{j=1}^{m}g(b_j,b_{-n})}
\frac{\prod_{j=1}^{m} g(b_j,b_{m+1})}
{\prod_{j=1}^{m}f(a_j,b_{m+1})}.
\end{align*}
Reformulate it by the convention. It leads to
\begin{align*}
\sum^m_{k=-n} f(a_k,b_k)\frac{\prod_{j=1}^{k-1}f(a_j,b_{-n})}
{\prod_{j=1}^{k}g(b_j,b_{-n})}
\frac{\prod_{j=m+1}^{k-1}g(b_j,b_{m+1})} {\prod_{j=m+1}^{k}
f(a_j,b_{m+1})} =-\frac{\prod_{j=1}^{m}f(a_j,b_{-n})}
{\prod_{j=1}^{m+1}g(b_j,b_{-n})} .
\end{align*}
Replace $m$ by $m-1$. It gives
\begin{align*}
\sum^{m-1}_{k=-n} f(a_k,b_k)\frac{\prod_{j=1}^{k-1}f(a_j,b_{-n})}
{\prod_{j=1}^{k}g(b_j,b_{-n})} \frac{\prod_{j=m}^{k-1}g(b_j,b_{m})}
{\prod_{j=m}^{k} f(a_j,b_{m})}
=-\frac{\prod_{j=1}^{m-1}f(a_j,b_{-n})}
{\prod_{j=1}^{m}g(b_j,b_{-n})},
\end{align*}
which is equivalent to
\begin{align*}
\sum^{m}_{k=-n} f(a_k,b_k) \frac{\prod_{j=m}^{k-1}g(b_j,b_{m})}
{\prod_{j=m}^{k} f(a_j,b_{m})}\frac{\prod_{j=1}^{k-1}f(a_j,b_{-n})}
{\prod_{j=1}^{k}g(b_j,b_{-n})}=0.
\end{align*}
As desired.

The next theorem gives the bilateral form of the $(f,g)$-inversion.
For convenience, write
$$
\frac{\prod_{j=1}^{\infty}f(a_j,c_j)} {\prod_{j=1}^{\infty}
f(a_j,d_j)}\frac{\prod_{j=1}^{m} g(b_j,d_j)} {\prod_{j=1}^{\infty}
g(b_j,c_j)}\qquad\mbox{for}\quad \lim_{m\mapsto
+\infty}\frac{\prod_{j=1}^{m}f(a_j,c_j)}
{\prod_{j=1}^{m}f(a_j,d_j)}\frac{\prod_{j=1}^{m} g(b_j,d_j)}
{\prod_{j=1}^{m}g(b_j,c_j)}.
$$
\begin{dl}\label{666} Let $f\bot g$, $M,N\in Z$,  and $\{A_M\}$ be an arbitrary sequence such that
 for all $A_M,A_N$, $$
\frac{\prod_{j=1}^{\infty}f(a_j,A_M)} {\prod_{j=1}^{\infty}
f(a_j,A_N)} \frac{\prod_{j=1}^{\infty} g(b_j,A_N)}
{\prod_{j=1}^{\infty}
g(b_j,A_M)}=\frac{\prod^{0}_{j=-\infty}f(a_j,A_N)}
{\prod^{0}_{j=-\infty}f(a_j,A_M)} \frac{\prod^{0}_{j=-\infty}
g(b_j,A_M)} {\prod^{0}_{j=-\infty}g(b_j,A_N)}
$$
and the limit
$$
h(M)=\lim_{N\mapsto M}\frac{\frac{\prod_{j=1}^{\infty}f(a_j,A_M)}
{\prod_{j=1}^{\infty} f(a_j,A_N)} \frac{\prod_{j=1}^{\infty}
g(b_j,A_N)} {\prod_{j=1}^{\infty}
g(b_j,A_M)}-\frac{\prod^{0}_{j=-\infty}f(a_j,A_N)}
{\prod^{0}_{j=-\infty}f(a_j,A_M)} \frac{\prod^{0}_{j=-\infty}
g(b_j,A_M)} {\prod^{0}_{j=-\infty}g(b_j,A_N)}}{g(A_M,A_N)}\neq 0.
$$
Then  a pair of matrices $F=(F_{n,k})_{n,k\in Z}$ and
$G=(G_{n,k})_{n,k\in Z}$  with entries given by
\begin{eqnarray*}
F_{n,k}&=&\frac{f(a_k,b_k)}{h(n)}\frac{\prod_{j=1}^{k-1}f(a_j,A_n)}
{\prod_{j=1}^{k}g(b_j,A_n)}\qquad\mbox{and}\qquad G_{n,k}=
\frac{\prod_{j=1}^{n-1}g(b_j,A_k)} {\prod_{j=1}^{n}
f(a_j,A_k)}\end{eqnarray*}
  is a
matrix inversion.
\end{dl}
{\sl Proof.} In (\ref{news111555}), set $m,n\mapsto \infty$. Then
\begin{align*}
\sum^\infty_{k=-\infty}
f(a_k,b_k)g(c_k,d_k)\frac{\prod_{j=1}^{k-1}f(a_j,c_j)}
{\prod_{j=1}^{k} f(a_j,d_j)}\frac{\prod_{j=1}^{k-1}g(b_j,d_j)}
{\prod_{j=1}^{k}g(b_j,c_j)}\nonumber\\=
\frac{\prod_{j=1}^{\infty}f(a_j,c_j)} {\prod_{j=1}^{\infty}
f(a_j,d_j)}\frac{\prod_{j=1}^{m} g(b_j,d_j)} {\prod_{j=1}^{\infty}
g(b_j,c_j)}- \frac{\prod^{0}_{j=-\infty}f(a_j,d_j)}
{\prod^{0}_{j=-\infty}f(a_j,c_j)} \frac{\prod^{0}_{j=-\infty}
g(b_j,c_j)} {\prod^{0}_{j=-\infty}g(b_j,d_j)}.
\end{align*}
Given two integers $M,N$, set further $c_j\mapsto A_M, d_j\mapsto
A_N$. Then
\begin{align}
\sum^\infty_{k=-\infty}
f(a_k,b_k)g(A_M,A_N)\frac{\prod_{j=1}^{k-1}f(a_j,A_M)}
{\prod_{j=1}^{k} f(a_j,A_N)}\frac{\prod_{j=1}^{k-1}g(b_j,A_N)}
{\prod_{j=1}^{k}g(b_j,A_M)}\nonumber\\=
\frac{\prod_{j=1}^{\infty}f(a_j,A_M)} {\prod_{j=1}^{\infty}
f(a_j,A_N)}\frac{\prod_{j=1}^{m} g(b_j,A_N)} {\prod_{j=1}^{\infty}
g(b_j,A_M)}- \frac{\prod^{0}_{j=-\infty}f(a_j,A_N)}
{\prod^{0}_{j=-\infty}f(a_j,A_M)} \frac{\prod^{0}_{j=-\infty}
g(b_j,A_M)} {\prod^{0}_{j=-\infty}g(b_j,A_N)}. \label{bilateral}
\end{align}
Define
\begin{align*}
h_{M,N}=\frac{1}{g(A_M,A_N)}
\left\{\frac{\prod_{j=1}^{\infty}f(a_j,A_M)} {\prod_{j=1}^{\infty}
f(a_j,A_N)} \frac{\prod_{j=1}^{\infty} g(b_j,A_N)}
{\prod_{j=1}^{\infty} g(b_j,A_M)}\right.\nonumber
\\\left.-
\frac{\prod^{0}_{j=-\infty}f(a_j,A_N)}
{\prod^{0}_{j=-\infty}f(a_j,A_M)} \frac{\prod^{0}_{j=-\infty}
g(b_j,A_M)} {\prod^{0}_{j=-\infty}g(b_j,A_N)}\right\}. \end{align*}
 Reformulate (\ref{bilateral}) with this note to get
\begin{align}
\sum^\infty_{k=-\infty} f(a_k,b_k)
\frac{\prod_{j=1}^{k-1}f(a_j,A_M)}{\prod_{j=1}^{k}g(b_j,A_M)}
\frac{\prod_{j=1}^{k-1}g(b_j,A_N)}{\prod_{j=1}^{k}
f(a_j,A_N)}=h_{M,N}.
\end{align}
The relation suggests that
\begin{align}
PG=H,\label{mine}
\end{align}
 where $P=(P_{M,k})$, $G=(G_{k,N})$, and $H=(h_{M,N})$ are three matrices with
 entries given by
\begin{eqnarray*}
P_{M,k}=f(a_k,b_k)
\frac{\prod_{j=1}^{k-1}f(a_j,A_M)}{\prod_{j=1}^{k}g(b_j,A_M)}\qquad\mbox{and}\quad
G_{k,N}= \frac{\prod_{j=1}^{k-1}g(b_j,A_N)} {\prod_{j=1}^{k}
f(a_j,A_N)},\end{eqnarray*} respectively. It should be possible that
$H$ is diagonal. Indeed, such $H$ exists under the known conditions,
whose entries are given by
\begin{eqnarray*}
h_{M,N}=\delta_{M,N}h(M),
\end{eqnarray*}
where
$$
h(M)=\lim_{N\mapsto M}\frac{\frac{\prod_{j=1}^{\infty}f(a_j,A_M)}
{\prod_{j=1}^{\infty} f(a_j,A_N)} \frac{\prod_{j=1}^{\infty}
g(b_j,A_N)} {\prod_{j=1}^{\infty}
g(b_j,A_M)}-\frac{\prod^{0}_{j=-\infty}f(a_j,A_N)}
{\prod^{0}_{j=-\infty}f(a_j,A_M)} \frac{\prod^{0}_{j=-\infty}
g(b_j,A_M)} {\prod^{0}_{j=-\infty}g(b_j,A_N)}}{g(A_M,A_N)}\neq 0.
$$
Thus, we get
$$
H^{-1}=\left(\frac
 {1}{h(M)}\delta_{M,N}\right).
$$
Finally, it follows from (\ref{mine}) that
$$
G^{-1}=H^{-1}P.
$$
 Equate the $(n,k)$-entries on both sides of this matric identity.
The result is proved.

As an important case is the following bilateral matrix inversion due
to Schlosser \cite{schloss}. He obtained this result from an
instance of Bailey¡¯s very-well-poised ${}_{6}\phi _{6}$ summation
theorem and used it successful to derive a lot of summation formulas
of bilateral hypergeometric series.
\begin{lz}
Let $a,b$ and $c$ be indeterminates. Then two infinite matrices
$A=(A_{n,k})$ and $B=(B_{n,k})$ are inverses of each other where
\begin{align}
A_{n,k}=\frac{(aq/b, bq/a, aq/c, cq/a, bq, q/b, cq, q/c;
q)_{\infty}} {(q, q,aq, q/a, aq/bc, bcq/a, cq/b, bq/c; q)_{\infty}}\nonumber\\
 \times\frac{ (1- bcq^{2n}/a)}{(1-
bc/a)}\frac{(b; q)_{n+k}(a/c;q)_{k-n}} {(cq; q)_{n+k} (aq/b;
q)_{k-n}}\end{align} and
\begin{align}
 B_{n,k}=\frac{(1-aq^{2n}) }{(1-
a)}\frac{(c;q)_{n+k}(a/b;q)_{n-k}}{(bq; q)_{n+k} (aq/c; q)_{n-k}}
q^{n-k}.\end{align}
\end{lz}
{\sl Proof.} It suffices to set in Theorem \ref{666}
$$
f(x,y)=g(x,y)=(y-x)(1-\frac{a}{bc}xy)
$$
and make the substitution
$$
A_n\mapsto q^{-n}, a_j\mapsto bq^j, b_j\mapsto cq^{j}.
$$
Some direct calculation gives
\begin{eqnarray*}
&&\prod_{j=-\infty}^{\infty}\frac{f(a_j,A_M)}{g(b_j,A_M)}=
\frac{(a/c,cq/a,b,q/b;q)_{\infty}}{(a/b,bq/a,c,q/c;q)_{\infty}};\\
&&\prod_{j=-\infty}^{\infty}\frac{f(b_j,A_N)}{g(a_j,A_N)}=
\frac{(a/b,bq/a,c,q/c;q)_{\infty}}{(a/c,cq/a,b,q/b;q)_{\infty}},
\end{eqnarray*}
 from which it follows that
$$
\frac{\prod_{j=1}^{\infty}f(a_j,A_M)} {\prod_{j=1}^{\infty}
f(a_j,A_N)} \frac{\prod_{j=1}^{\infty} g(b_j,A_N)}
{\prod_{j=1}^{\infty}
g(b_j,A_M)}=\frac{\prod^{0}_{j=-\infty}f(a_j,A_N)}
{\prod^{0}_{j=-\infty}f(a_j,A_M)} \frac{\prod^{0}_{j=-\infty}
g(b_j,A_M)} {\prod^{0}_{j=-\infty}g(b_j,A_N)}
$$
and
\begin{eqnarray*}
h(M)&=&(q, q,aq, q/a, aq/bc,
bcq/a,cq/b,bq/c;q)_{\infty}\\
&&\frac{q^{3M}(1-
bc/a)(1-a/(cq^M))(1-a/(bq^M))}{(1-b)(1-a)(1-a/b)(1-a/c)(1-
bcq^{2M}/a)}.
 \end{eqnarray*}
Further calculation leads to
\begin{eqnarray*}
&&P_{n,k}=\frac{(c-b)(1-aq^{2k})(cq;q)_{n}
(a/(cq);q^{-1})_{n-1}}{(b;q)_{n+1}(a/b;q^{-1})_{n}}\frac{(b;
q)_{n+k}(a/c;q)_{k-n}} {(cq; q)_{n+k} (aq/b;
q)_{k-n}}q^{n+k},\\
&&G_{n,k}=\frac{q^k(a/(bq);q^{-1})_{k-1}(bq;q)_{k}}
{(a/c;q^{-1})_k(c;q)_{k+1}}\frac{(c;q)_{n+k}(a/b;q)_{n-k}}{(bq;
q)_{n+k} (aq/c; q)_{n-k}}.
\end{eqnarray*}
Observe that there exist the following relations
$$
P=LAT,\quad G=T^{-1}BS,
$$
where $L, T$, and $S$ are three diagonal matrices with the
corresponding diagonal entries given by
\begin{eqnarray*}
L_{k,k}&=&\frac{q^k(c-b)(cq;q)_{k}
(a/(cq);q^{-1})_{k-1}}{(b;q)_{k+1}(a/b;q^{-1})_{k}}\frac{(1- bc/a)}{
(1- bcq^{2k}/a)}\\
&&\frac{(q, q,aq, q/a, aq/bc, bcq/a, cq/b, bq/c; q)_{\infty}}{(aq/b,
bq/a, aq/c, cq/a, bq, q/b, cq, q/c;
q)_{\infty}};\\
 T_{k,k}&=&q^k(1-aq^{2k});\\
S_{k,k}&=&\frac{q^{2k}(a/(bq);q^{-1})_{k-1}(bq;q)_{k}}
{(1-a)(a/c;q^{-1})_k(c;q)_{k+1}},\,\qquad\mbox{respectively}.
\end{eqnarray*}
Since that $PG=H$ gives
$$
LABS=H,
$$
direct calculation states that $H=LS$, so $AB=(\delta_{n,k})$. The
desired result follows.

 Note that the calculation for $h(M)$ provides an alternative proof of the famous transformation (cf. \cite[Ex. 5.21]{10})
\begin{eqnarray*}
&&(aq/b, bq/a, aq/c, cq/a, bq, q/b, cq, q/c; q)_{\infty}(1/b, 1/c, 1/d, bcd/a^2; q)_{\infty}\\
&-&(aq/b, aq/c, aq/d, bcdq/a; q)_{\infty}(b/a, c/a, d/a, a/bcd; q)_{\infty}\\
&=&(aq, q/a, aq/bc, bcq/a, aq/bd, bdq/a, cdq/a, aq/cd; q)_{\infty}.
\end{eqnarray*}

\section*{Acknowledgements}
The author wants to thank Prof.W.C.Chu for his  providing the
references \cite{chu,chu1,chu2,chu3}. Thanks are also due to the
anonymous referee for his very detailed comments and Prof.L.Zhu  for
helpful suggestions.

  \bibliographystyle{amsplain}

\end{document}